\documentclass[onecolumn,10pt,cleanhead,cleanfoot]{asme2e}

\usepackage{graphicx}
\usepackage{amsmath}
\usepackage{amssymb}
\usepackage{multirow}
\usepackage{xcolor}
\usepackage{graphicx} 
\usepackage{setspace}
\usepackage{mathrsfs}

\title{Convergence Acceleration of Favre-Averaged Non-Linear Harmonic Method}

\author{Feng Wang\thanks{Corresponding author, email: feng.wang@eng.ox.ac.uk}
    \affiliation{
	Oxford Thermo-Fluids Institute\\
	Department of Engineering Science\\
	University Of Oxford \\
	UK
    }
}

\author{Kurt Weber
    \affiliation{
	Compressor Aero Method\\
	Rolls-Royce Corporation\\
	USA
    }
}

\author{David Radford
    \affiliation{
	Noise Engineering\\
	Rolls-Royce plc\\
	UK
    }
}

\author{Luca di Mare
    \affiliation{
	Oxford Thermo-Fluids Institute\\
	Department of Engineering Science\\
	University Of Oxford \\
	UK
    }
}

\author{Marcus Meyer
    \affiliation{
	CFD Methods\\
	Rolls-Royce Deutschland Ltd \& Co KG.\\
	Germany
    }
}

\begin{document}

\maketitle    

\begin{abstract}
{ \it This paper develops a numerical procedure to accelerate the convergence of the Favre-averaged Non-Linear Harmonic (FNLH) method. The scheme provides a unified mathematical framework for solving the sparse linear systems formed by the mean flow and the time-linearized harmonic flows of FNLH in an explicit or implicit fashion. The approach explores the similarity of the sparse linear systems of FNLH and leads to a memory efficient procedure, so that its memory consumption does not depend on the number of harmonics to compute. The proposed method has been implemented in the industrial CFD solver Hydra. Three test cases are used to conduct a comparative study of explicit and implicit schemes in terms of convergence, computational efficiency, and memory consumption.  Comparisons show that the implicit scheme yields better convergence than the explicit scheme and is also roughly 7 to 10 times more computationally efficient than the explicit scheme with 4 levels of multigrid. Furthermore, the implicit scheme consumes only approximately $50\%$ of the memory required by the explicit scheme with four levels of multigrid. Compared with the full annulus unsteady Reynolds averaged Navier-Stokes (URANS) simulations, the implicit scheme produces comparable results to URANS with computational time and memory consumption that are two orders of magnitude smaller.   }

\end{abstract}

\section*{Introduction}

Computational Fluid Dynamics (CFD) has been routinely used in the aerodynamic design of fluid machinery. A fluid machinery (i.e. turbomachinery) can involve a series of rotational and stationary rows of blades to exchange energy between the machine and the working fluid. Although the exchange of energy between the working fluid and the machine is dependent on the unsteadiness of the flow~\cite{comp_aero_book}, steady Reynolds Averaged Navier-Stokes (RANS) has been the workhorse in the industry to design turbomachinery. However, unsteady interactions among turbomachinery components can be significant and should be considered at the design stage. Unsteady RANS (URANS) simulation is an approach to model these unsteady interactions, but they can increase the computational cost by several orders of magnitude compared to the steady RANS approach. The industry has shown significant interest in developing computationally efficient design methods to model these unsteady interactions. This work is a step in this direction. 
           
If the effect of stochastic flow perturbations on the mean flow can be approximated by a suitable turbulence model, the unsteady flow in a turbomachine can be effectively decomposed into a time-mean flow and an unsteady periodic perturbation.  Periodic fluctuations are deterministic and can be computed efficiently by harmonic methods, which transform unsteady time domain problems into steady problems in the frequency domain. This removes the constraint of global time stepping in the nonlinear time marching approach (i.e. URANS). For multistage simulations, phase-lag boundaries can be applied to the periodic boundaries, and this effectively avoids the need to model multiple passages or the whole wheel in the simulations. Nevertheless, this approach ignores the interactions between stochastic and periodic perturbations; besides, the frequencies of the periodic perturbations normally need to be known a priori. 
        
This work concerns the family of harmonic methods that couples a mean flow solver and a time-linearized NS flow solver in the frequency domain. The coupling between the mean flow and the perturbations can be evaluated using the flow perturbations of the linearized NS solver and the mean flow variables.  He and Ning~\cite{He1998} proposed the first type of such approach, which is called the Non-Linear Harmonic (NLH) method, but the authors did not elucidate the conversion among time-averaged conservative and primitive variables in the mean flow.  Wang and di Mare~\cite{fnlh} proposed a more elaborate and rigorous approach to include the nonlinear coupling between the mean flow and the flow perturbations for compressible flows, and the approach is termed the Favre-averaged NLH (FNLH) method. FNLH can be used to simulate unsteady rotor-stator interactions~\cite{fnlh}, multi-row~\cite{fnlh_multirow} interactions and more recently multi-shaft~\cite{Wang2023} unsteady interactions. Furthermore, FNLH can also be used to efficiently simulate non-uniform flow fields in an assembly of bladerows with non-uniform geometries~\cite{sfnlh_aiaa,fan_sfnlh_jot}. The other family of harmonic methods is the Harmonic Balance (HB) method~\cite{Hall_hb_2002}. This approach does not require a linearized NS solver, but it solves a finite set of steady-state solutions at fixed time levels in the time domain simultaneously. All of these time-level solutions are coupled through periodic boundaries and a spectral source term. Readers can refer to Hall et al.~\cite{Hall2013} for more details on this approach. 

Harmonic methods (i.e., FNLH) use pseudo-time stepping to solve the unsteady problem in the Fourier space~\cite{fnlh}. For each pseudo-time step, the mean flow and the time-linearized flow need to be computed. Convergence acceleration techniques are crucial to the computational performance of harmonic methods. NLH used an explicit scheme with geometric multigrid to accelerate its convergence~\cite{He1998}. However, there is no previous research in the public domain on developing an implicit scheme for NLH or FNLH. On the other hand, it is worth mentioning that there has been previous work on developing implicit schemes to accelerate the convergence of the HB method, the fully implicit approach~\cite{Woodgate_2009,Su_2009} has been proposed but memory usage can increase rapidly as more harmonics are computed. Considering the fact that the coupling among the harmonics and the mean flow can be weak, a segmented approach~\cite{Sicot2008,Frey2014} can be used to solve the mean flow and each harmonic perturbation separately, so that the memory requirement depends weakly on the number of harmonics to compute. However, since the HB method does not require a dedicated time-linearized solver, these approaches cannot be used directly to accelerate the convergence of FNLH or NLH. 

There has been extensive research on the acceleration of convergence of a steady NS flow solver ~\cite{JAMESON1991,Yoon1988,saad_book} and also a time-linearized flow solver~\cite{Campobasso2003} using explicit or implicit approaches. This work will bridge the gap between convergence acceleration techniques in these two aspects and develop an implicit scheme to accelerate the convergence of FNLH. This is the first contribution of this paper. The second contribution of this paper is to formulate a unified formulation to solve the sparse linear systems of FNLH in an explicit or implicit fashion. The third contribution is to conduct a comparative study of explicit and implicit schemes with respect to convergence, computational efficiency, and memory usage. The proposed method is implemented in the industrial CFD solver Hydra~\cite{moinier_thesis}, which a major aero-engine manufacturer routinely uses for aerothermal design.

The paper is organized as follows: the numerical algorithms of Hydra are briefly introduced in the first place, and the FNLH method is then briefly described. The generalized RK scheme for FNLH is then detailed. The performance of this generalized explicit/implicit RK scheme is demonstrated. This is followed by conclusions and future work. 


\section*{Brief Description of Hydra}

\begin{figure}
\centering 
\includegraphics[height=6cm]{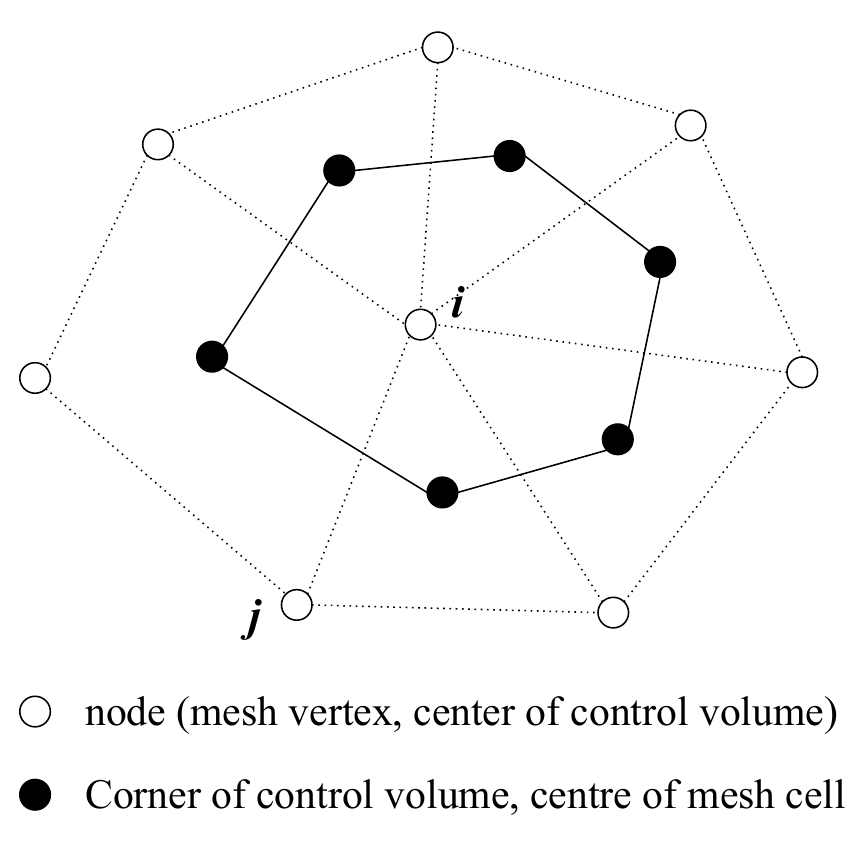}
\caption{Control volume of an internal node in a node-based finite volume scheme.}
\label{fig:51_control_volume}       
\end{figure}

The numerics of Hydra are briefly described here, and readers can refer to Ref.~\cite{Hydra_mg,Campobasso2003,moinier_thesis} for more details. Hydra is a finite-volume CFD solver and uses the edge-based data structure, and this can be illustrated by Fig.~\ref{fig:51_control_volume}.  The solver is second-order accurate in space. Inviscid flux is calculated with the Jameson-Schmidt-Turkel (JST) scheme using matrix dissipation~\cite{JAMESON1981}. The viscous flux is computed using the central difference. To calculate the flow gradient, the gradient on each edge $(i,j)$ is first calculated using a mixture of the central difference and the average flow gradient at nodes $i$ and $j$. Gradient estimates for each node are given by the Green-Gauss gradients~\cite{moinier_thesis}. In Hydra, a set of turbulence models is implemented. In this work, the Spalart-Allmaras (SA)~\cite{spalart1992one} turbulence model with wall function is used. With respect to time integration,  A low-dissipation and low-dispersion RK scheme RK($s, p$) is implemented in Hydra, where $s$ is the number of stages and $p$ is the number of evaluations of the viscous contribution per time step. In Hydra $s=5$ and $p=3$. This scheme has been widely used in the literature and has demonstrated good performance~\cite{Swanson2007,moinier_thesis,Turkel1993} . To accelerate the convergence of the nonlinear solver, the geometric multigrid~\cite{Moinier2002} and the incomplete matrix factorization~\cite{misev_thesis} have been implemented. The former leads to an explicit RK scheme, and the latter leads to an implicit RK scheme. 

Hydra has an existing time-linearized harmonic flow solver. TAPENADE~\cite{Hascoet2013} is used to perform the source code transformation to automatically differentiate the nonlinear solver and produce the linearized code. Only the subroutines that contribute to the non-linear residuals are differentiated, but the time-stepping for the time-linearized solver is still coded manually and similar to the non-linear solver the RK(5,3) scheme is used. The linearized code can only work with real numbers; therefore, the real and imaginary parts of harmonic flow perturbations are stored separately.  With respect to turbulence modeling, although source code transformation can also be applied to subroutines that are related to turbulence models, frozen turbulence is still used in this work for enhanced robustness. Future work is planned to remove this restriction.

In the baseline Hydra, the time-linearized harmonic flow solver is independent from the non-linear solver. Furthermore, it can only compute one harmonic at a time for a single bladerow. As FNLH requires one to compute the mean flow and the harmonic flow simultaneously for multiple harmonics and multiple bladerows,  the limitations of the baseline Hydra have been completely removed in the current code development to successfully implement the FNLH method.   



\section*{Favre-averaged Non-Linear Harmonic Method}
The Favre-averaged Non-Linear Harmonic (FNLH) method is a computational framework that efficiently couples the mean flow and its finite-amplitude periodic perturbations for compressible flows in the frequency domain.  Here we provide a brief description of FNLH and the details of FNLH can be found in Wang and di Mare~\cite{fnlh,fnlh_multirow,sfnlh_aiaa}. 

For a flow that is subject to periodic disturbance, the primitive variables of the $i^{\text{th}}$ control volume $\mathbf{W}_i=(\rho, u, v, w, p)^T$ can be decomposed as:
\begin{equation}
    \mathbf{W}_i = \hat{\mathbf{W}}_i + {\mathbf{W}}''_i
    \label{eqn:flow_decomp}
\end{equation}
in which $\hat{\mathbf{W}}_i$ is the Favre-averaged primitive variable and ${\mathbf{W}}''_i$ is the flow perturbation relative to the Favre-averaged value. It is noted that, for Favre-averaging, density $\rho$ and pressure $p$ are still time-averaged~\cite{wilcox_book}. Their perturbations are denoted as $\rho'$ and $p'$. The conversion between Favre-averaged and time-averaged conservative variables $\bar{\mathbf{Q}}_i=(\bar{\rho}, \overline{\rho u}, \overline{\rho v}, \overline{\rho w}, \overline{\rho E})^T$ is trivial~\cite{cebeci_book}, for example $\overline{\rho u} = \bar{\rho} \tilde{u}$.

\begin{figure}
\centering 
\includegraphics[height=5cm]{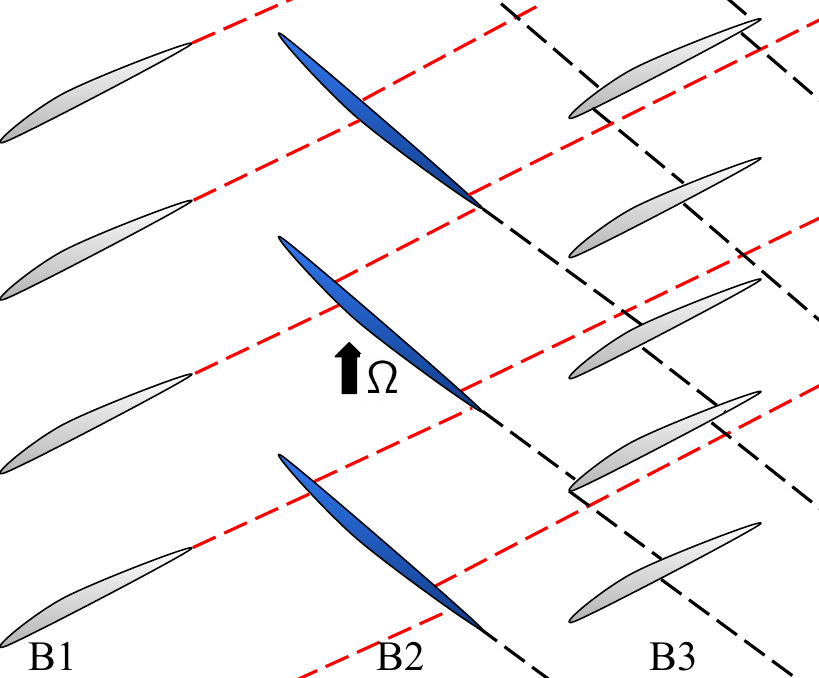}
\caption{Schematic of multistage flows in turbomachinery.}
\label{fig:wake_diag_new}       
\end{figure}

In the context of periodic unsteady flows in turbomachinery, considering the schematic in Fig.~\ref{fig:wake_diag_new}, which features a stator(B1) - rotor(B2) - stator(B3) setting. The flow field of B3 can be decomposed to the combination of a passage-averaged time-mean value, an unsteady disturbance, and a stationary passage-to-passage disturbance. If we consider the stationary disturbance as a "slowly moving" unsteady disturbance and its frequency is approaching zero. The flow field in B3 can be decomposed as a passage-averaged time-mean value $\tilde{\phi}$ and an unsteady flow disturbance $\phi''$, whose frequency can be zero. Here, $\phi$ is a dummy flow variable.

If the perturbation of a dummy variable ${\phi}''_i$ of a control volume can be represented as the harmonic form:
\begin{equation}
\label{eqn::phi_fft}
\phi''_i = \sum^{l=N_h}_{l=1} (\hat{\phi}_{i,l}e^{I \omega_{i,l} t} + \hat{\phi}_{i,-l}e^{-I \omega_{i,l} t})
\end{equation}
where $l$ is the harmonic index, $N_h$ is the number of harmonics, $I$ is the imaginary unit $\sqrt{-1}$, $\omega_{l}$ is the angular frequency of the $l^{\text{th}}$ harmonic and $t$ is time. Using the flow decomposition and the representation of the flow perturbation in the harmonic form, the coupled system of steady mean flow and harmonic flow perturbations in the frequency domain for a control volume can be obtained and solved with pseudo-time stepping. This can be represented symbolically as:
\begin{align}
   \mathbf{M}\frac{\partial \mathbf{W}_i}{\partial t^*} + \mathbf{R}_i + \text{DF}_i &= 0  \label{eqn::mean_system}  \\ 
  \mathbf{M}\frac{\partial \mathbf{\hat{W}}_{l,i}}{\partial t^*} + \mathbf{\hat{R}}_{l,i} &= 0 \label{eqn::lin_system} 
\end{align}
in which $t^*$ is the pseudo-time, $l$ is the harmonic index, $V$ is the volume of the current control volume and $\mathbf{M}$ is the Jacobian of the transformation from primitive variable $\mathbf{W}$ to conservative variables $\mathbf{Q}$. $\mathbf{R}$ is the non-linear residual of the mean flow and $\mathbf{\hat{R}}_{l}$ is the non-linear residual for the $l^{th}$ harmonic of the time-linearized flow. DF represents the contribution of deterministic fluxes~\cite{df_jot}, which takes into account the coupling between mean flow and harmonic flow perturbations. and can be determined using the mean flow variable and harmonic flow perturbations~\cite{fnlh}.  The nonlinear residual of the time-linearized flow $\mathbf{\hat{R}}_{l,i}$ of control volume $i$ is:
\begin{equation}
    \label{eqn::lin_residual}
    \mathbf{\hat{R}}_{l,i} = I \omega_{l,i} \mathbf{M} V_i \mathbf{\hat{W}}_{l,i}  + \frac{\partial\mathbf{R}}{\partial \mathbf{W}} \mathbf{\hat{W}}_{l,i}
\end{equation}
in which $\frac{\partial\mathbf{R}}{\partial \mathbf{W}} \mathbf{W}''$ represents the linearized flux that is accumulated on all the faces of the control volume of node $i$, here the indices (i.e. $j$) are ignored and it is written symbolically for simplicity. In Hydra, the linearized flux is computed by performing an automatic differentiation of the subroutines that are related to the nonlinear residuals.

It should be noted that, although Hydra solves the conservative form of the NS equation, it works with primitive variables. This is why the coupled system in Equation~\ref{eqn::mean_system} and Equation~\ref{eqn::lin_system} is written in the form of primitive variables. The formulation of FNLH is preferred in Hydra over the NLH method, because the NLH method is formulated with time-averaged conservative variables and requires an extra step to convert time-averaged conservative variables to time-averaged primitive variables~\cite{He1998}, but this has been avoided in FNLH and will lead to significantly less modifications of the baseline code. 

\begin{figure}
\centering 
\includegraphics[height=7cm]{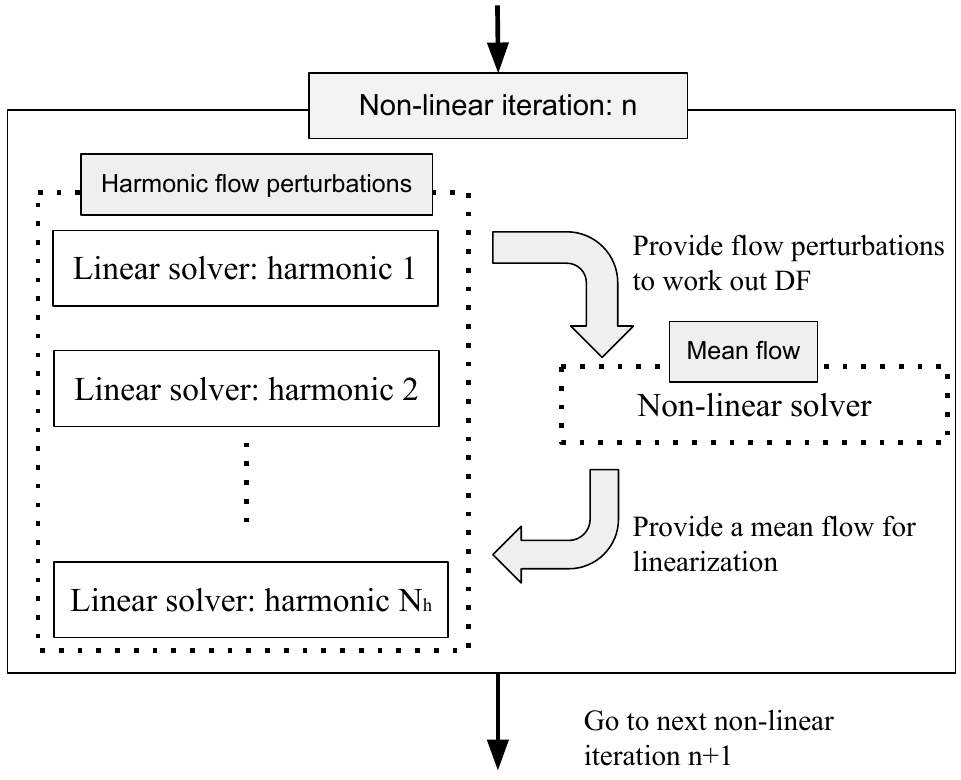}
\caption{Flow chart of FNLH for each non-linear iteration.}
\label{fig:fnlh_flowchart}       
\end{figure}

Similarly to the baseline steady flow solver, FNLH uses pseudo-time stepping to solve the mean flow (Equation~\ref{eqn::mean_system}) and the time-linearized harmonic flow (Equation~\ref{eqn::lin_system}). Each pseudo-time step is called a non-linear iteration or outer iteration. In this non-linear iteration, the harmonic flow perturbations are evaluated in the first place, DF is then computed using the available flow perturbations in the harmonic form, and the mean flow is then updated with the inclusion of DF. This procedure is illustrated by Fig.~\ref{fig:fnlh_flowchart}. It is noted that the cross-coupling of the harmonics is ignored in the current work. This will be explored in our future work.

In a multistage setting of FNLH, bladerow interfaces for the mean flow and the time-linearized harmonic flow are treated separately. For the mean flow, the original mixing plane formulation~\cite{Denton1992} is corrected by DF to take into account the unsteady effect~\cite{fnlh}. For the time-linearized flow, a unified formulation has been developed to transmit disturbances in a multi-row environment~\cite{fnlh_multirow}.  Quasi-3D non-reflective treatment~\cite{giles_nrbc_report} are applied to the mean flow and the time-linearized flow.


\section*{Convergence Acceleration Procedure of FNLH with The Generalized RK scheme}
Hydra uses the Runge-Kutta scheme to advance the solution in pseudo-time. For the baseline steady flow solver,  a generalized Runge-Kutta (RK) scheme has been developed for both explicit~\cite{Moinier2002} and implicit~\cite{Misev2018} schemes, but this procedure is not available to the time-linearized harmonic flow solver. In this paper such a generalized RK scheme is extended for FNLH to solve the mean flow and the time-linearized harmonic flows together. This will lead to a unified formulation to accelerate the convergence of FNLH. In the following, the generalized RK scheme for the mean-flow solver is described in the first place, and this is then extended to solve the time-linearized harmonic flow.




\subsection*{Generalized RK Scheme for The Mean Flow}
The generalized RK scheme in Hydra is based on the hybrid five-stage scheme RK (5,3)~\cite{Swanson2007}, where the non-linear viscous residuals are evaluated at stages 1, 3 and 5. It can be written as:

\begin{align}
\label{eqn:rk}
  \mathbf{Q}^{n,(0)} & = \mathbf{Q}^n \nonumber \\ 
  \mathbf{Q}^{n,(k)} & = \mathbf{Q}^{n,(0)} - \Gamma \alpha^{(k)} \mathbf{K}^{-1}  (\mathbf{R}^{(k)} - \mathbf{R}^*), k = \{1,2\cdots5\} \nonumber \\
  \mathbf{Q}^{n+1} & = \mathbf{Q}^{n,(5)} 
\end{align}

In Equation~\ref{eqn:rk}, without losing generality, the control volume index $i$ has been ignored for simplicity, $n$ is the current pseudo-time step index, $k$ is the RK stage index, $\mathbf{R}^{(k)}$ is the non-linear residual of the $k^{\text{th}}$ RK stage, $\mathbf{R}^*$ is the additional term related to multigrid operation, the formulation of which can be found in Moinier~\cite{moinier_thesis}. For each RK stage, a solution update needs to be computed to advance the solution. Such a solution update can be worked out in an explicit way, and this is the classic explicit RK approach. In addition, one can also invert the Jacobian matrix in an approximate way to obtain the update of the solution for each RK stage. In this case, since there is no multigrid operation involved, $\mathbf{R}^*$ in Equation~\ref{eqn:rk} is set to zero. The matrix inversion can be realized with the incomplete lower-upper decomposition with zero level of fill-in. This technique can be used as a preconditioner for Newton-Krylov methods, and it can also be used as a smoother to work out solution updates.  The generalized RK scheme allows both explicit and implicit schemes to be written in a unified way, according to Misev~\cite{misev_thesis} the formulations of $\Gamma$ and $\mathbf{K}$ for explicit and implicit approaches can be written as:
\begin{equation}
    \Gamma = \left\{
\begin{array}{ll}
       \sigma_{\text{CFL}} & \text{Explicit} \\
      \frac{1}{\varepsilon} & \text{Implicit}
\end{array} 
\right.
\end{equation}
 
\begin{equation}
\label{eqn:precon}
    \mathbf{K} = \left\{
\begin{array}{ll}
      \mathbf{P} & \text{Explicit} \\
      \frac{\mathbf{P}}{\varepsilon \sigma_{\text{CFL}}} + \mathbf{J} & \text{Implicit}
\end{array} 
\right.
\end{equation}
where $\mathbf{J}$ is the residual Jacobian matrix $\frac{\partial{\mathbf{R}}}{\partial \mathbf{Q}}$. The residual Jacobian is assembled using the information of the immediate neighbors with first-order spatial accuracy, and hence it is an approximation of the true Jacobian.  $\mathbf{P}$ is the diagonal block of the Jacobian matrix $\mathbf{J}$. For the explicit scheme,  $\Gamma$ is the the Courant–Friedrichs–Levy (CFL) number $\sigma_{CFL}$, in Hydra $\sigma_{CFL}=2$ is recommended and used by default in this work. It should be noted that for the implicit scheme, $\Gamma$ no longer represents the CFL number, but it is related to a relaxation factor $\varepsilon$. The detailed derivation can be found in Ref.~\cite{Misev2018,Swanson2007}, a value of 0.6 is suggested in the previous work~\cite{Misev2018}.

To improve the convergence rate and stability, the non-linear residual of the $k^{\text{th}}$ stage $\mathbf{R}^{(k)}$ is decomposed into two components, namely the inviscid contribution $\mathbf{R}^{(k),inv}_i$ and the dissipative or viscous contribution $\mathbf{R}^{(k),vis}$:
\begin{align}
  \mathbf{R}^{(k)} &= \mathbf{R}^{(k),inv} + \mathbf{R}^{(k),vis} 
\end{align}
in which, the inviscid and viscous contributions can be evaluated as:
\begin{align}
  \mathbf{R}^{(k),inv} &= \mathbf{R}^{inv} (\mathbf{Q}^{n,(k-1)})  \\
  \mathbf{R}^{(k),vis} &= \beta_k \mathbf{R}^{vis}(\mathbf{Q}^{n,(k-1)}) + (1-\beta_k) \mathbf{R}^{(k-1),vis} 
\end{align}

$\mathbf{R}^{inv} (\mathbf{Q}^{n,(k-1)})$ is the inviscid contribution using the flow variable $\mathbf{Q}^{n,(k-1)}$, $\mathbf{R}^{(k),vis}$ is the viscous contribution and is a mixture of the viscous flux using the flow variable $\mathbf{Q}^{n,(k-1)}_i$ and the viscous residual $\mathbf{R}^{(k-1),vis}$ from the previous RK stage. The coefficients $\alpha^{(k)}$ and $\beta^{(k)}$ are summarized in Table~\ref{tab:rk_coff}. These coefficients are optimized to damp high-frequency errors and lead to better convergence when used with multigrid~\cite{swanson_turkel_rk_coeff}. 

\begin{table}
\centering
\caption{Coefficients of the hybrid RK(5,3) scheme}
\label{tab:rk_coff}

\begin{tabular}{ c | ccccc }
\hline \hline
 k & 1 & 2 & 3 & 4 & 5 \\ 
 \hline
 $\alpha^{(k)}$ & $\frac{1}{4}$ & $\frac{1}{6}$ & $\frac{3}{8}$ & $\frac{1}{2}$ & 1  \\ \hline 
 $\beta^{(k)}$ & 1 & 0 & $\frac{14}{25}$ & 0 & $\frac{11}{25}$ \\   
\hline \hline
\end{tabular}
\end{table}

\subsection*{Generalized RK Scheme for Time-Linearized Harmonic Flows}

With respect to the time-linearized flow, the structure of the generalized RK scheme for the mean flow can be re-used. $\mathbf{Q}^{n,(k)}$ should be replaced with the harmonic flow perturbations $\mathbf{\hat{Q}}^{n,(k)}_l$, and $\mathbf{R}^{(k)}$ should be replaced with the non-linear residuals from the time-linearized flows $\mathbf{\hat{R}}^{(k)}_l$, for example:

\begin{align}
\label{eqn:rk_lin}
  \mathbf{\hat{Q}}^{n,(0)}_l & = \mathbf{\hat{Q}}^n_l \nonumber \\ 
  \mathbf{\hat{Q}}^{n,(k)}_l & = \mathbf{\hat{Q}}^{n,(0)} - \Gamma \alpha^{(k)} \mathbf{\hat{K}}^{-1}_l  (\mathbf{\hat{R}}^{(k)}_l - \mathbf{\hat{R}}^*_l), k = [1:5] \nonumber \\
  \mathbf{\hat{Q}}^{n+1}_l & = \mathbf{\hat{Q}}^{n,(5)}_l 
\end{align}
in which $\Gamma$ is the same as the mean flow RK scheme, some attention is required for the matrix preconditioner. Due to the spectral source term in the time-linearized solver, a diagonal matrix $I \omega_{l} \mathbf{I} V$ is added to the diagonal block of $\mathbf{K}$, where $\omega_l$ is the angular frequency. This leads to the preconditioner matrix for each harmonic $\mathbf{\hat{K}}_l$, which consists of complex numbers, and this is shown in Equation~\ref{eqn:precon_linear}. 

\begin{equation}
    \label{eqn:precon_linear}
    \mathbf{\hat{K}}_l = \left\{
\begin{array}{ll}
      \mathbf{P} + I \omega_l \mathbf{I} V   & \text{Explicit} \\
      \frac{\mathbf{P}}{\varepsilon \sigma_{\text{CFL}}} + \mathbf{J} + I \omega_l \mathbf{I} V  & \text{Implicit}
\end{array} 
\right.
\end{equation}
in which $I$ is the imaginary unit, $V$ is the volume of a control volume, and $\mathbf{I}$ is the unit matrix whose size is $5 \times 5$ for a 3D problem. Because $\omega_l$ and the interblade phase angle are generally different for each harmonic, $\mathbf{\hat{K}}_l$ is then different from harmonic to harmonic, but the sparse pattern of $\mathbf{\hat{K}}_l$ remains the same as the mean flow. This feature can be explored to develop a memory-efficient approach to solving FNLH.


\subsection*{Solution Update}

For each RK stage of the mean flow, an update of the solution $\Delta \mathbf{Q}^{(k)}$ to the vector of the intermediary solution $\mathbf{Q}^{n,(k)}$ must be calculated, and this reads as follows:
\begin{equation}
\label{eqn:rk_update}
    \Delta \mathbf{Q}^{(k)} = \mathbf{Q}^{n,(k)} - \mathbf{Q}^{n,(0)} = -\Gamma \alpha^{(k)} \mathbf{K}^{-1} (\mathbf{R}^{(k)}_i - \mathbf{R}^*)
\end{equation}

To use more traditional symbols for the discussion of solving sparse linear systems, Equation~\ref{eqn:rk_update} can be written in a more familiar form $\mathbf{A} \mathbf{x} = \mathbf{B}$. In the following, $\mathbf{A}$ refers to the left-hand side (LHS) of the sparse linear system and $\mathbf{B}$ refers to the right-hand side of this system. For the mean flow, they represent:
\begin{align}
  \mathbf{A} &=  (\Gamma \alpha^{(k)}\mathbf{K}^{-1})^{-1} = \frac{\mathbf{K}}{\Gamma \alpha^k} \label{eqn::lhs_mean} \\
  \mathbf{B} &= -(\mathbf{R}^{(k)} - \mathbf{R}^*) \\
  \mathbf{x} &= \Delta \mathbf{Q}^{(k)}
\end{align}
For the explicit scheme, $\mathbf{A}$ is a block diagonal matrix, and its inversion can be performed point by point. This is essentially a block Jacobi method, and inverting the block-diagonal matrix is a simple and fast process.  For the implicit scheme, $\mathbf{A}$ is a large sparse matrix. Solving the linear system exactly is seldom economical, because the improvement of the convergence rate is not good enough to justify the high computational cost of inverting this matrix exactly. In this work, the procedure similar to that proposed by Swanson et al.~\cite{Swanson2007} is followed. This matrix is inexactly inverted, and the linear system is solved iteratively. The incomplete lower-upper decomposition with zero level of fill-in (which is commonly denoted as ILU(0)) is used to approximate the matrix inversion. To be more specific, if an ILU(0) procedure is applied to the matrix $\mathbf{A}$. i.e. $\mathbf{A} \approx \mathbf{L} \mathbf{U}$. $\mathbf{L}$ is a lower triangular matrix and $\mathbf{U}$ is an upper triangular matrix. By decomposing the original matrix to a lower and upper triangular matrix, the matrix inversion can be computed efficiently. The zero level of fill-in means that the sparsity pattern of $\mathbf{L}$ and $\mathbf{U}$ is chosen to be the same as the sparsity pattern of the matrix $\mathbf{A}$, which means if an entry in $\mathbf{A}$ is zero, then that entry in the matrix of $\mathbf{L}$ or $\mathbf{U}$ is also zero. Therefore, $\mathbf{L}\mathbf{U}$ is just an approximation of $\mathbf{A}$, and so will the matrix inversion by ILU(0). Because there is no need to solve the linear system exactly, ILU(0) is a popular choice for matrix inversion to solve the linear system.


With respect to the time-linearized harmonic flow, after inspecting Equation~\ref{eqn:precon_linear} and Equation~\ref{eqn:precon}, it can be seen that the LHS of the linear system of the linearized flow has the same sparse pattern as the one for the nonlinear flow. The difference comes mainly from the diagonal block. Furthermore, if there are periodic boundaries in the simulation, phase-lagged boundary conditions are applied on the periodic boundaries. Since the interblade phase angle will be different for each harmonic, this will also lead to a difference of the LHS for each harmonic, but it has no impact on the sparse pattern of the LHS. The linear system for the $l^{th}$ harmonic can be written as:
\begin{equation}
    \frac{\mathbf{K}_l}{\Gamma \alpha^{k}} \hat{\mathbf{x}}_l = \hat{\mathbf{B}}_l 
\end{equation}
in which $\mathbf{\hat{B}}_l$ and $\mathbf{\hat{x}}_l$ denote the RHS and solution update of the linear system of the $l^{\text{th}}$ harmonic, respectively. Furthermore, the LHS of this system can be written using that of the mean flow as:
\begin{equation}
    \label{eqn::lhs_lin}
    (\mathbf{A} + \frac{I \omega_l \mathbf{I} V} {\Gamma \alpha^k}) \hat{\mathbf{x}}_l = \hat{\mathbf{B}}_l 
\end{equation}

When comparing the LHS of the linear system of the mean flow (Equation~\ref{eqn::lhs_mean}) and that of the linearized flow(Equation~\ref{eqn::lhs_lin}), it can be seen that they share the same sparse pattern. This means that the residual Jacobian matrix $\mathbf{J}$ only needs to be assembled once by the nonlinear flow for each nonlinear iteration, and is then reused by the linearized flow. From an implementation point of view, the mean flow and the first harmonic of the linearized flow will allocate memory for its LHS, respectively, but the LHS of all the harmonics of the linearized flow reuses the memory space of that of the first harmonic. This treatment is used for both explicit and implicit schemes and allows the memory usage of FNLH to not depend on the number of harmonics to compute; therefore, the resulting procedure is memory efficient. It is worth mentioning that a similar procedure was also explored in the harmonic balance method~\cite{Frey2014}, however, this procedure is described here in the context of FNLH.

\subsection*{Complex Number Representation}
The sparse linear system of the time-linearized flow consists of complex numbers. We can cast this complex linear system (i.e. Equation~\ref{eqn::lhs_lin}) into its equivalent system with real numbers, or we can directly solve this complex linear system. 
\begin{equation}
\label{eqn:matrix_z_real}
\begin{bmatrix}
    \mathbf{A}  &  -{\omega}_l \mathbf{I}\mathbf{V} /(\Gamma \alpha^k)     \\
    {\omega}_l \mathbf{I}\mathbf{V}/(\Gamma \alpha^k)  &  \mathbf{A}      
\end{bmatrix}
\begin{bmatrix}
    \text{Re}\{\Delta \hat{\mathbf{x}}_l\}   \\
    \text{Im}\{\Delta \hat{\mathbf{x}}_l\}       
\end{bmatrix}
= 
\begin{bmatrix}
    \text{Re} \{ \hat{\mathbf{B}}_l\}      \\
    \text{Im} \{ \hat{\mathbf{B}}_l\}     
\end{bmatrix} 
\end{equation}

The benefit of using real numbers is that the existing implementation of incomplete matrix factorization can be reused. However, as shown in Equation~\ref{eqn:matrix_z_real}, the sparse pattern of the new LHS differs from that of the mean flow solver and this means that a new lookup table is required for the matrix factorization (i.e. ILU(0)). For the approach using complex numbers, the sparse pattern of the LHS is the same as that of the mean flow solver, and this also facilitates parallel implementation (i.e. parallel domain decomposition) of incomplete matrix factorization. However, a complex-number version of the baseline ILU(0) procedure needs to be implemented. From the authors' point of view, turning an incomplete matrix factorization into its complex version is much easier than implementing a new look-up table and its related domain decomposition. Therefore, in this work, the approach with complex numbers is chosen.    




\subsection*{Parallel Implementation}
\begin{figure}
\centering 
\includegraphics[height=8cm]{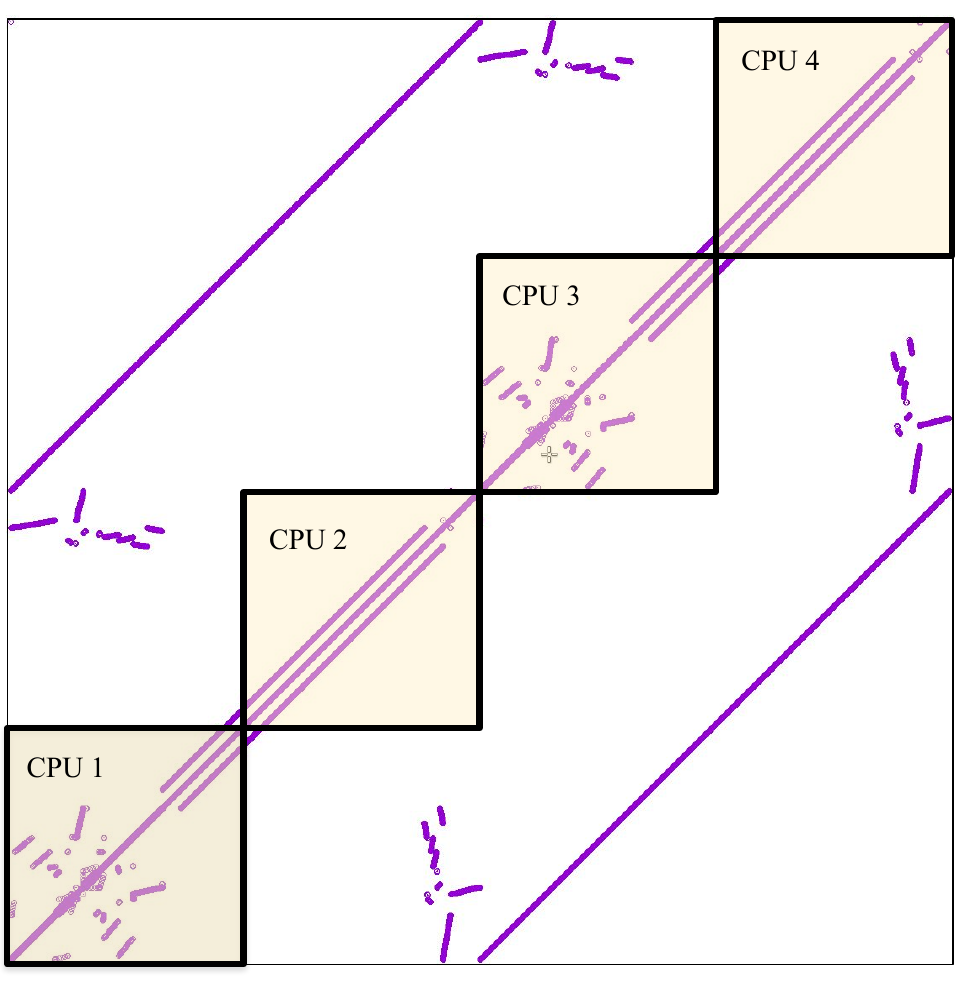}
\caption{Parallel domain decomposition for incomplete matrix factorization.}
\label{fig:partition_matrix}       
\end{figure}

For parallel implementation of the implicit scheme for both the mean flow and the harmonic flows, a domain decomposition strategy is required to partition the LHS of the sparse linear system of the mean flow (Equation~\ref{eqn::lhs_mean}) or the linearized flow (Equation~\ref{eqn::lhs_lin}) for different processors. As the sparse patterns of the linearized flow and the mean flow are identical, the same domain decomposition strategy is used.   In this work, the partition is done in such a way that the entries related to two nodes that are in different partitions are ignored. This can be demonstrated in Fig.~\ref{fig:partition_matrix}. It shows the sparse pattern of the LHS for the mean flow of a compressor blade passage that is meshed with a multi-block structured grid, and its partition of the LHS on four processors. Each of the 4 boxes encloses a sub-linear system for each processor.  These sub-linear systems only couple nodes that are partitioned into the same processor. This effectively leads to a parallelisable set of sub-problems for each processor. 

In addition, as the coupling of the nodes that are in different processors is ignored, the convergence rate of the implicit solver will degrade as more processors are used. For the extreme case, where the number of processors is equal to the number of nodes, the parallel implicit solver will degrade to the explicit scheme.  It is found that serious convergence degradation can occur if the total number of nodes to the number of processors is less than $O(10^3)$~\cite{Misev2018}. On the other hand, to obtain a good scalability for parallel computation, this ratio should normally be at least $O(10^4)$. Therefore, in practice, the convergence rate will only have a weak dependence on the number of processors, and this effect will be demonstrated in the test cases.   

\section*{Results}
We demonstrate the computational performance of the explicit and implicit versions of the generalized RK scheme for FNLH and perform a comparative study of both approaches in terms of convergence, wall clock time and memory consumption. Three test cases are used to serve this purpose. The first case is to model the rotor-rotor interactions in a compressor, the second case is to trace hot streak migrations in turbine stages, and the third test case is the low-pressure compression system of a turbofan engine.  



The key parameters of the generalized RK scheme are summarized in Table~\ref{tab:fnlh_para}. For the implicit solver, the incomplete LU factorization with zero fill-in level is used and is denoted as ILU(0). A fixed threshold of 0.01 is used for the linear residual to iteratively solve the sparse linear system for both the mean flow and the time-linearized harmonic flow. Taking the mean flow as an example, 0.01 means that $||\mathbf{B}-\mathbf{A}\mathbf{x}||_2$ drops by two orders of magnitude compared to its initial value. In addition, the maximum number of iterations to solve the sparse linear systems is 10.  For the explicit scheme, a V cycle is used for the multigrid and the order of spatial accuracy for the coarse grid levels is set to 1. 

\begin{table}[]
\centering
\caption{Key Parameters of the generalized RK scheme}
\label{tab:fnlh_para}
\begin{tabular}{lll}
\hline\hline
Scheme & Component & Value  \\
\hline
Implicit & ILU fill levels  &  0    \\
         & Linear residual convergence tolerance  & 1e-2      \\
         & Linear iteration limit  &10      \\
Explicit & MG Cycle  &  V-cycle    \\
         & Spatial order on coarse grid  & 1      \\
\hline \hline
\end{tabular}
\end{table}

\subsection*{Compressor Rotor-Rotor Interaction}
This test case is the mid-span stream tube of the front 1.5 stages of an 8-stage high-speed machine~\cite{simod_jot}, and the stream tube accounts for approximately $4.8\%$ of the span height. The first- and second-stage rotors of this 8-stage compressor are used~\cite{fnlh_multirow}. Blade counts are scaled to the ratio of 3(R1):5(S1):4(R2) to reduce the computational cost of URANS. The in-house meshing tool~\cite{mesh_jpp} is used to create multiblock structured meshes to discretize the computational domain in the blade-to-blade section and will be used to mesh the computational domain for all test cases. In the spanwise direction, one layer of hexahedral mesh is used, and the total hexahedral element is 33090. The average $y^{+}$ on the blade surface is around 5. The Sparlart-Allmaras turbulence model~\cite{SPALART1992} is used by default in all test cases. Non-reflective treatment for mean flow~\cite{Saxer1993} and linearized flow~\cite{giles_nrbc_report} is applied at the bladerow interfaces. For the URANS simulation, single-passage meshes are replicated in the circumferential direction to form a $30^{\circ}$ sector. Dual-time stepping is used and 500 time steps are used for this simulation of the $30^{\circ}$ sector, which corresponds to 100 time steps for the stator to sweep one passage of the rotor that has the smaller pitch.

Table~\ref{tab::fnlh_setup_r1_s1_r2} shows the FNLH setup for this compressor case, and $\Omega$ is the shaft speed.  15 harmonics are used in total. In the current FNLH implementation in Hydra, all bladerows compute the same total number of harmonics. For R1, 15 harmonics are used to calculate the potential fields of S1. In the ``note" column, it shows that these modes are applied on the exit boundary of R1 and are used to calculate the S1 potential field. For S1, 12 harmonics are used to model the R1 wake and they are applied at the inlet of S1. 3 harmonics are used to compute the potential field of R2 and are applied at the exit boundary of S1. This information is also provided in the ``note " column.  For R2, 8 harmonics are used for the wake of S1 and 7 harmonics are used to model the wake of R1. For the harmonic representation of the R1 wake in R2, its wave number is a linear combination of the blade counts of R1 and S1, and its temporal frequency is related to the shaft speed $\Omega$. These harmonics are organized into three groups based on the value of $m$ and are referred to as ``wake-1", ``wake-2" and ``wake-3", respectively.   This FNLH setup was used in our previous work, and for more details, the readers can refer to the previous work~\cite{fnlh_multirow}.

\begin{table}
\caption{FNLH Setup for the compressor case}
\label{tab::fnlh_setup_r1_s1_r2}
\centering
\tabcolsep=0.11cm
\begin{tabular}{cccc} 
 \hline
 \hline
 Bladerow & Temporal mode & Spatial mode & Note \\ 
 \hline
 R1 & $i\text{N}_{S1}\Omega, i\in\{1,2,...15\}$ & $i\text{N}_{S1}, i\in\{1,2,...15\}$ & exit, S1 potential field \\ 
       \hline
 S1 & $i\text{N}_{R1}\Omega, i\in\{1,2,...12\}$ & $i\text{N}_{R1}, i\in\{1,2,...12\}$ & inlet, R1 wake \\ 
      & $-i\text{N}_{R2}\Omega, i\in\{1,2,3\}$ & $i\text{N}_{R2}, i\in\{1,2,3\}$ & exit, R2 potential field \\ 
\hline
 R2 & $i\text{N}_{S1}\Omega, i\in\{1,2,...8\}$ & $i\text{N}_{S2}, i\in\{1,2,...8\}$ & inlet, S1 wake \\ 
      & $-n\text{N}_{S1}\Omega, n \in \{0,1,2\}$ & $m\text{N}_{R1} + n \text{N}_{S1}, m \in \{1\}$ & inlet, R1 wake-1 \\ 
      & $-n\text{N}_{S1}\Omega, n \in \{0,1\}$ & $m\text{N}_{R1} + n \text{N}_{S1}, m \in \{ 2 \}$ & inlet, R1 wake-2 \\ 
      & $-n\text{N}_{S1}\Omega, n \in \{0,1\}$ & $m\text{N}_{R1} + n \text{N}_{S1}, m \in \{ 3 \}$ & inlet, R1 wake-3 \\ 
\hline 
\hline
\end{tabular}
\end{table}
 
\begin{figure}
\centering 
\includegraphics[height=6.5cm]{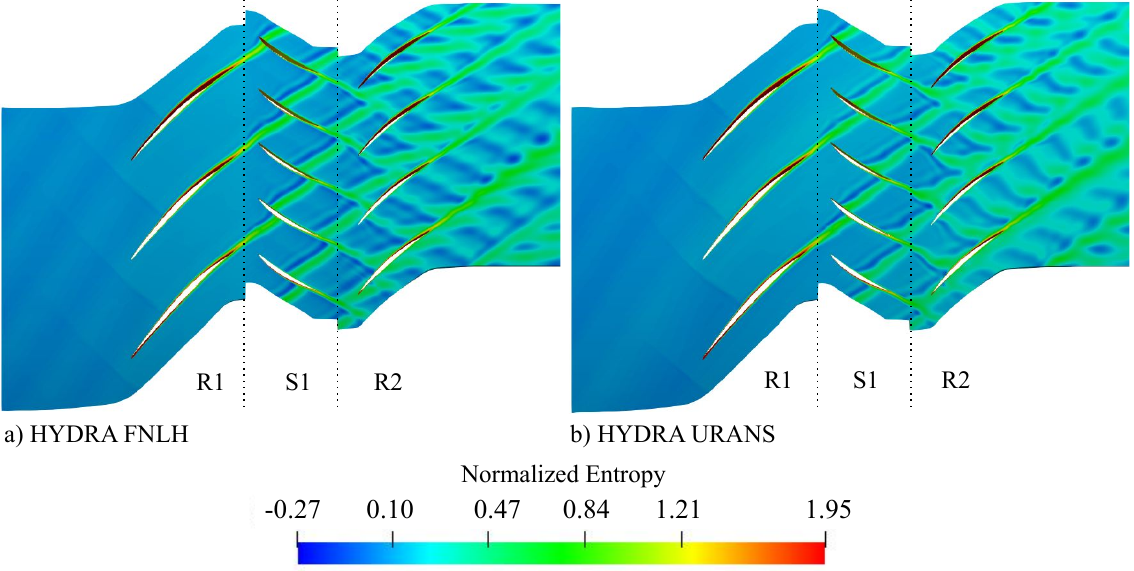}
\caption{Instantaneous entropy field of the rotor-stator-rotor case.}
\label{fig:r1_s1_r2_entropy}       
\end{figure}

\begin{figure}
\centering 
\includegraphics[height=9cm]{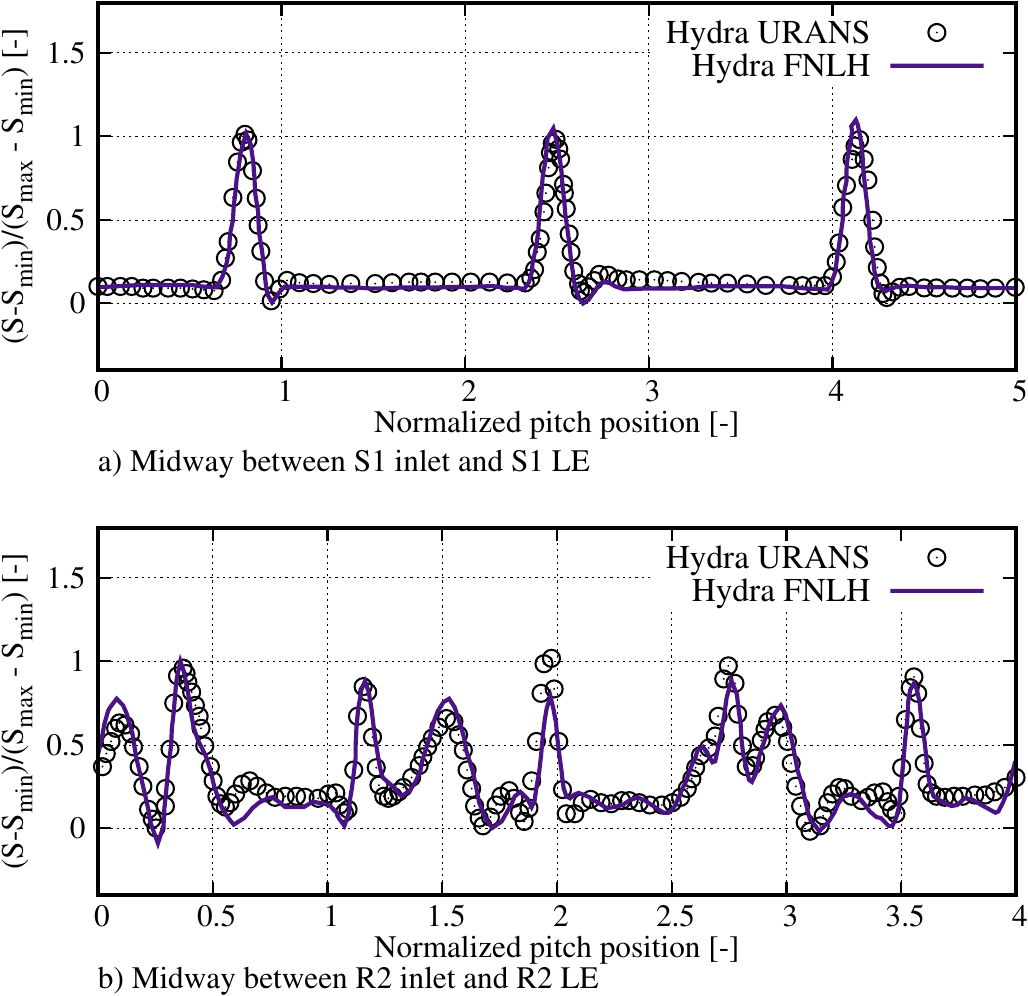}
\caption{Circumferential distribution of instantaneous entropy from FNLH.}
\label{fig:r1_s1_r2_s}       
\end{figure}

Figure~\ref{fig:r1_s1_r2_entropy} shows the instantaneous entropy fields for this 1.5-stage transonic compressor computed by FNLH and URANS. The entropy is normalized as ($S-S_{min})/(S_{max}-S_{min})$. $S_{min}$ and $S_{max}$ are the minimum and maximum entropy in the circumferential direction at a plane that is midway between the S1 inlet and the S1 LE.  FNLH shows good agreement with URANS, in particular, the wake of R1 is satisfactorily reconstructed in R2. Figure~\ref{fig:r1_s1_r2_s} shows the circumferential distributions of the normalized entropy at two streamwise locations, which are located midway between the S1 inlet and LE, the R2 inlet and LE, respectively. Good agreement with URANS data is observed. The results are computed by the implicit RK scheme, and hence we are confident that the generalized RK scheme still produces the correct results. Nevertheless, some discrepancies between FNLH and URANS data can still be observed; this is because frozen turbulence is assumed and cross-coupling is not included; future work is to address these issues and enhance the accuracy of the Hydra FNLH solver.     

\begin{figure}
\centering 
\includegraphics[height=7.0cm]{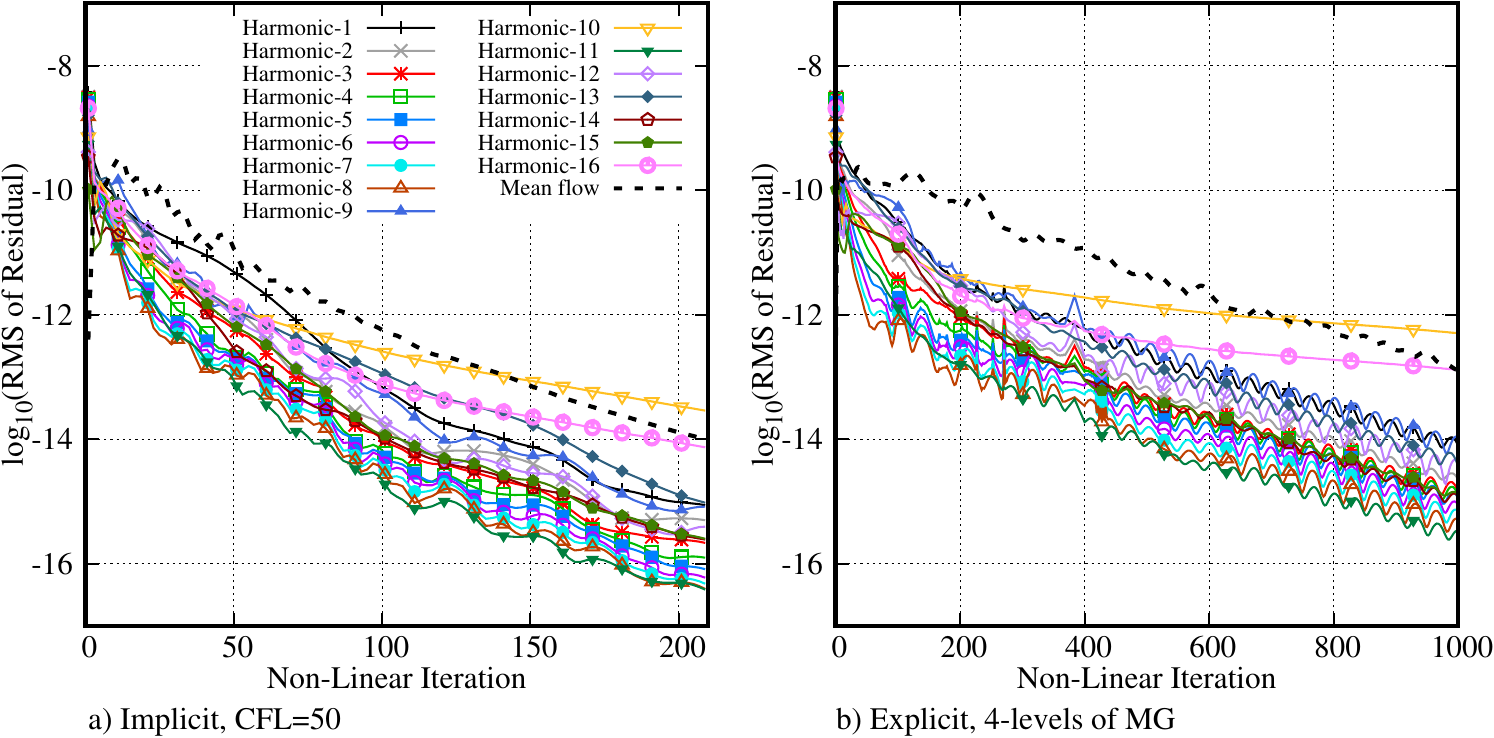}
\caption{Convergence history of all the harmonics with the implicit and explicit schemes for the 1.5-stage compressor case.}
\label{fig:r1_s1_r2_resd_all_h}       
\end{figure}

Figure~\ref{fig:r1_s1_r2_resd_all_h} shows the RMS of the residuals for the mean flow and the harmonic flows from the explicit and implicit schemes. For example, the residual for "Harmonic-0" is the residual of the first harmonic computed by \emph{all} the blade rows. According to the FNLH setup in Table~\ref{tab::fnlh_setup_r1_s1_r2}, Harmonic [0-7] represents the adjacent row interactions for all the blade rows. Harmonic [8-14] still computes the adjacent row interaction for R1 and S1, but for R2 they model the R1 wake in the computational domain of R2. Figure~\ref{fig:r1_s1_r2_resd_all_h} shows that for both implicit and explicit schemes, the harmonics converge at difference rates, and the harmonic that has a lower convergence rate in an explicit scheme also generally has a lower rate in the implicit scheme. In particular, for adjacent row interactions, harmonics with larger indices have shorter wave lengths and higher frequencies, and Fig.~\ref{fig:r1_s1_r2_resd_all_h} shows that these harmonics tend to have faster convergence. This could be explained by the fact that, because of the shorter wave length, the wave length to grid size ratio is smaller, and this means that these harmonics are more dissipated by the numerical scheme. 

Since a dozen harmonics could be computed in an FNLH run, it is not convenient to inspect the residuals for each harmonic. Therefore, the RMS of the residuals for all harmonics is monitored. This is defined as:
\begin{equation}
    R_{z} = \sqrt{\frac{1}{N_{h}} (R^2_{z,1} + R^2_{z,2} + R^2_{z,3} + \cdots R^2_{z,N_h})}
\end{equation}
where $R_{z,i}, i \in [1:N_h],$ is the RMS of the non-linear residual for the $i^{\text{th}}$  harmonic. In the following, the residual of the mean flow and the RMS of the residuals for all harmonics, namely $R_z$, is used to compare the performance of the explicit and implicit RK scheme.

\begin{figure}
\centering 
\includegraphics[height=6.5cm]{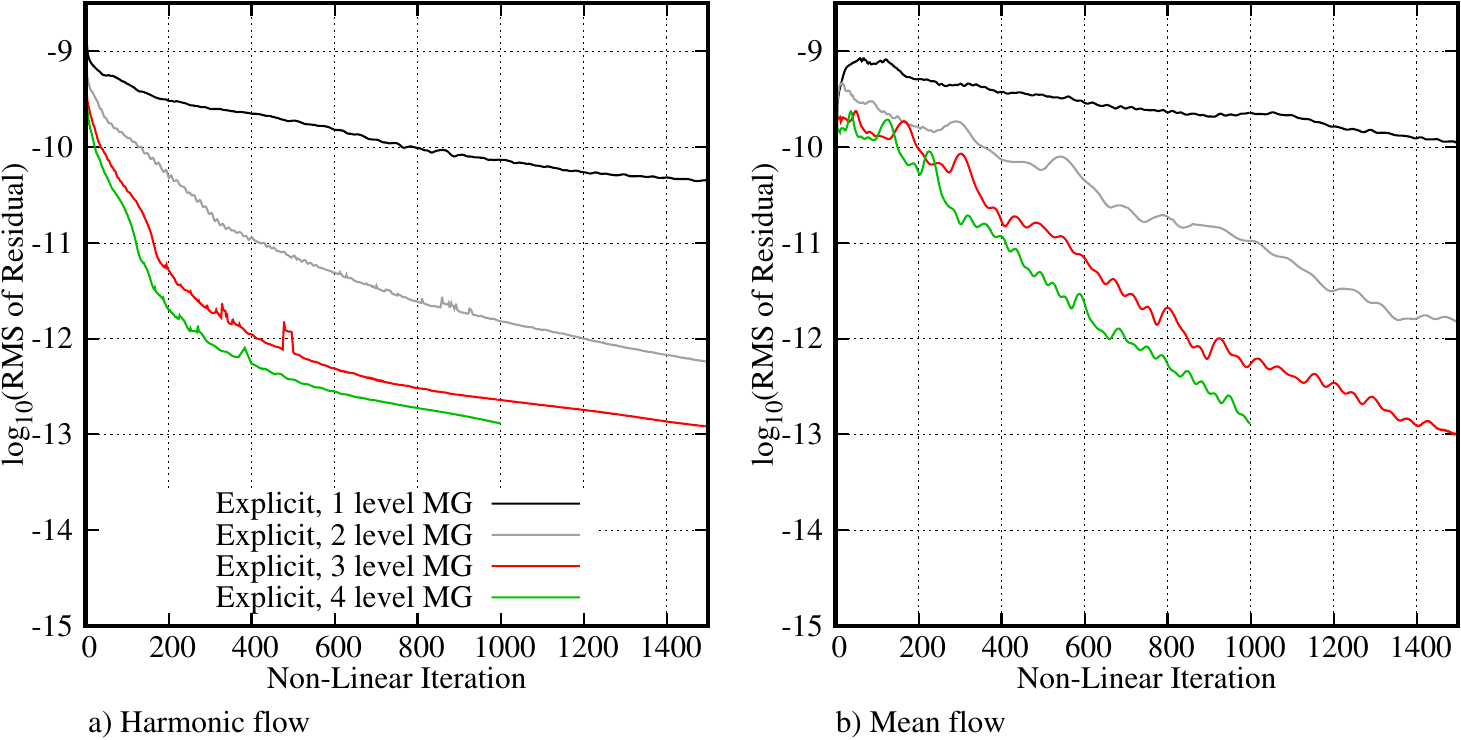}
\caption{Convergence history of explicit schemes with different levels of MG for the 1.5-stage compressor case.}
\label{fig:r1_s1_r2_resd_exp_mg}       
\end{figure}

For the performance of the explicit scheme with different levels of MG, Fig.~\ref{fig:r1_s1_r2_resd_exp_mg} shows the convergence history of FNLH for this compressor case. The explicit RK scheme with the V-cycle multigrid is used, and only the residual of the finest grid level is shown. It can be seen that using only one multigrid level yields poor convergence. As the number of MG levels increases, the convergence rate is significantly improved, but the benefit of adding more MG levels is also decreasing.  The non-linear mean flow seems to benefit more from multigrid compared to the time-linearized harmonic flow. The oscillations of residuals between 250 - 500 iterations are likely caused by the transonic flows in the rotors. As the residuals are further reduced and the flow becomes more established; such an oscillation is no longer present. 

It should be noted that as the number of MG levels increases, each nonlinear iteration becomes more expensive. It is also important to examine the residual history against wall-clock time. Because the performance of a CFD code will depend on the computer system, one way to overcome this problem is to employ a benchmark such as the TauBench code to nondimensionalize the computational cost (see Wang et al.~\cite{Wang2013}) and express the wall-clock time in terms of the Work Unit (WU), so that the running time of different CFD codes that run on different computers could be compared. In this work,  WU is computed as $\frac{N_{\text{proc}}T_{\text{FNLH}}}{T_1}$. $N_{\text{proc}}$ is the number of CPU cores that are used and $T_{\text{FNLH}}$ is the wall clock time of the FNLH simulation. $T_1$ is the wall clock time for executing the TauBench code in the local system with the following parameters: ``mpirun -np 1 ./TauBench -n 250000 -s 10", and these parameters are widely used in the CFD community.

\begin{figure}
\centering 
\includegraphics[height=6.5cm]{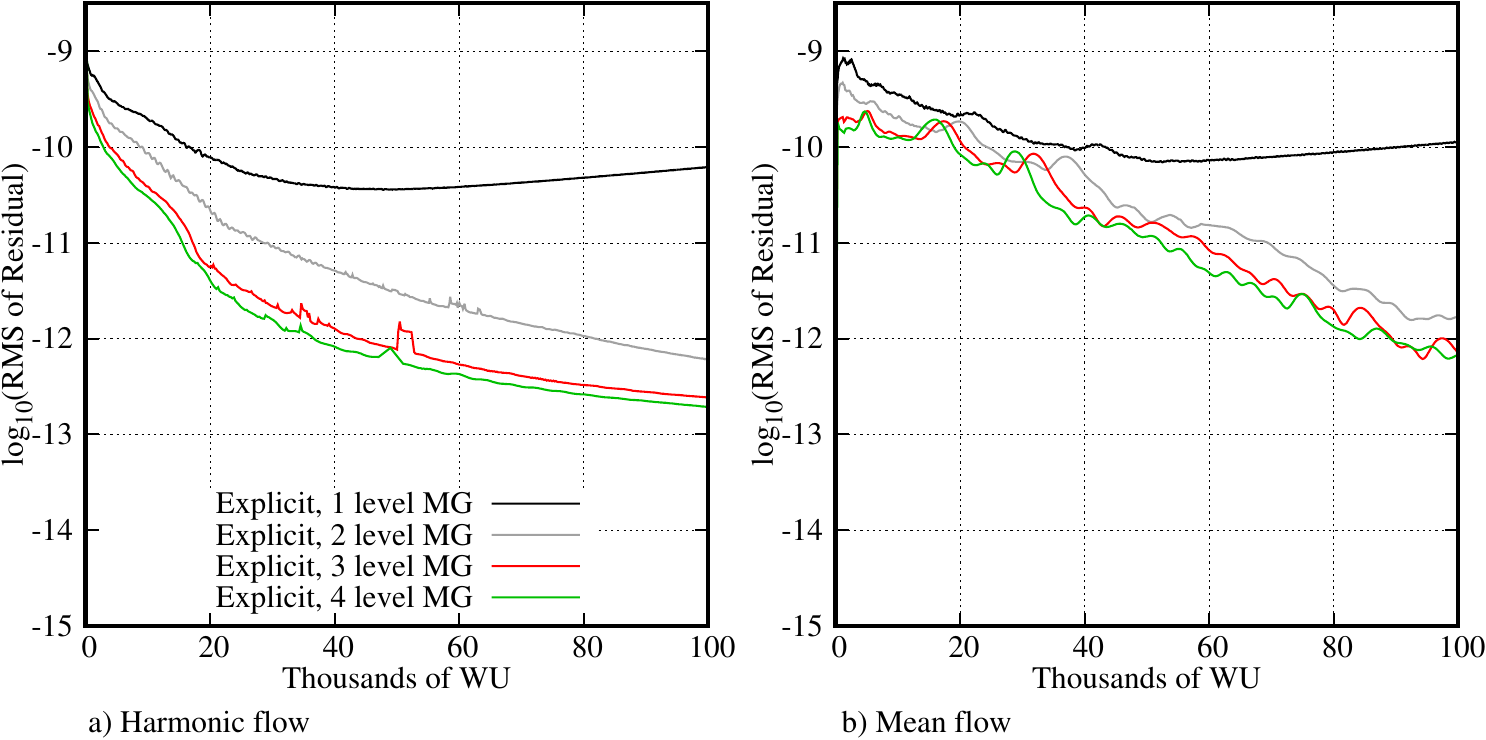}-\caption{Computational efficiency of explicit schemes with different levels of MG for the 1.5 stage compressor case.}
\label{fig:r1_s1_r2_resd_exp_mg_wu}       
\end{figure}

Figure~\ref{fig:r1_s1_r2_resd_exp_mg_wu} shows the WU against the nonlinear residual. For the explicit scheme, although each nonlinear iteration becomes expensive as the number of MG levels increases, 4 levels of MG is still the most computationally efficient one, and using only one level of MG is not computationally efficient and also shows poor convergence. In the following, the explicit scheme with 4 levels of MG will be used as the default explicit scheme and compared with the implicit RK schemes.

\begin{figure}
\centering 
\includegraphics[height=6.5cm]{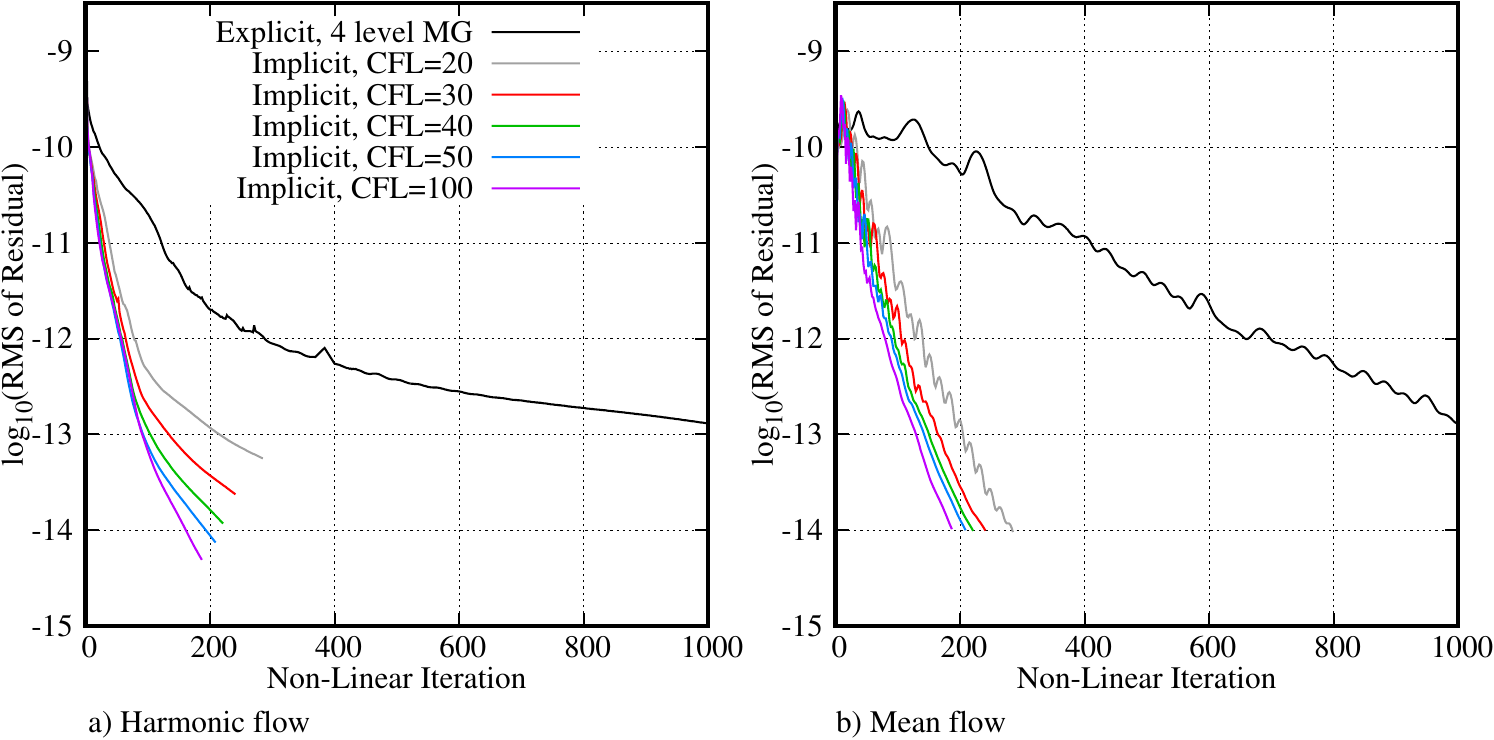}
\caption{Convergence history of implicit schemes with different CFL numbers for the 1.5-stage compressor case}
\label{fig:r1_s1_r2_resd_exp_imp}       
\end{figure}

Figure~\ref{fig:r1_s1_r2_resd_exp_imp} shows the comparison of the convergence histories of 4 levels of MG and implicit RK schemes with different CFL numbers. It can be seen that, in terms of non-linear iterations, implicit RK schemes show significant improvement over the explicit scheme with 4 levels of MG. In addition, the oscillations of the residuals observed between 200 and 500 iterations are not present in the implicit scheme. This is because the implicit scheme uses a full Jacobian matrix with first-order spatial accuracy for the preconditioner matrix shown in Equation~\ref{eqn:precon} and Equation~\ref{eqn:precon_linear} while the explicit scheme only uses the diagonal block of this matrix. Therefore, the solution update could be in a more accurate direction for the implicit scheme to reduce the residual. When there is strong non-linearity in the flow, this is beneficial to improve the convergence and also the robustness of the solution.

\begin{figure}
\centering 
\includegraphics[height=6.5cm]{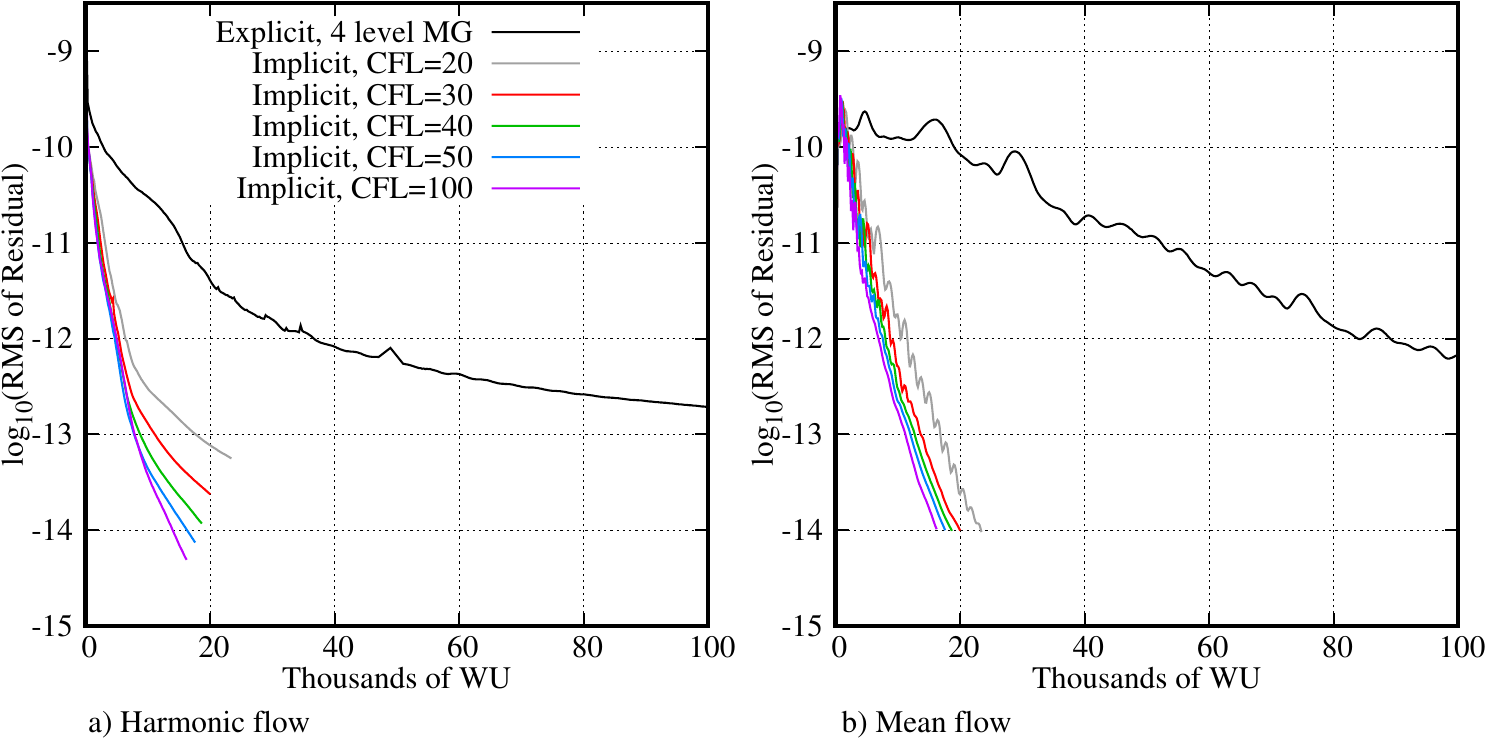}
\caption{Computational efficiency of implicit schemes for the 1.5-stage compressor case.}
\label{fig:r1_s1_r2_resd_exp_imp_wu}       
\end{figure}

In the implicit scheme, the CFL number acts as a relaxation factor. A smaller CFL number leads to better diagonal dominance, but this will also lead to slower convergence. Therefore, a judicious choice of the CFL number will strike a good balance in terms of convergence speed, stability, and computational efficiency. Figures~\ref{fig:r1_s1_r2_resd_exp_imp} and~\ref{fig:r1_s1_r2_resd_exp_imp_wu} show that as the CFL number gradually increases from 10 to 100, the convergence rate improves, but the benefit of increasing the CFL number also decreases rapidly. In addition, increasing the CFL reduces the diagonal dominance and could lead to potential robustness issues for complex flows in practice.  From our numerical experiment,  a CFL number of 40 or 50 is found to be a good choice, and this is consistent with the previous experience of developing the implicit scheme for the non-linear mean flow~\cite{misev_thesis}. 

\begin{figure}
\centering 
\includegraphics[height=6.5cm]{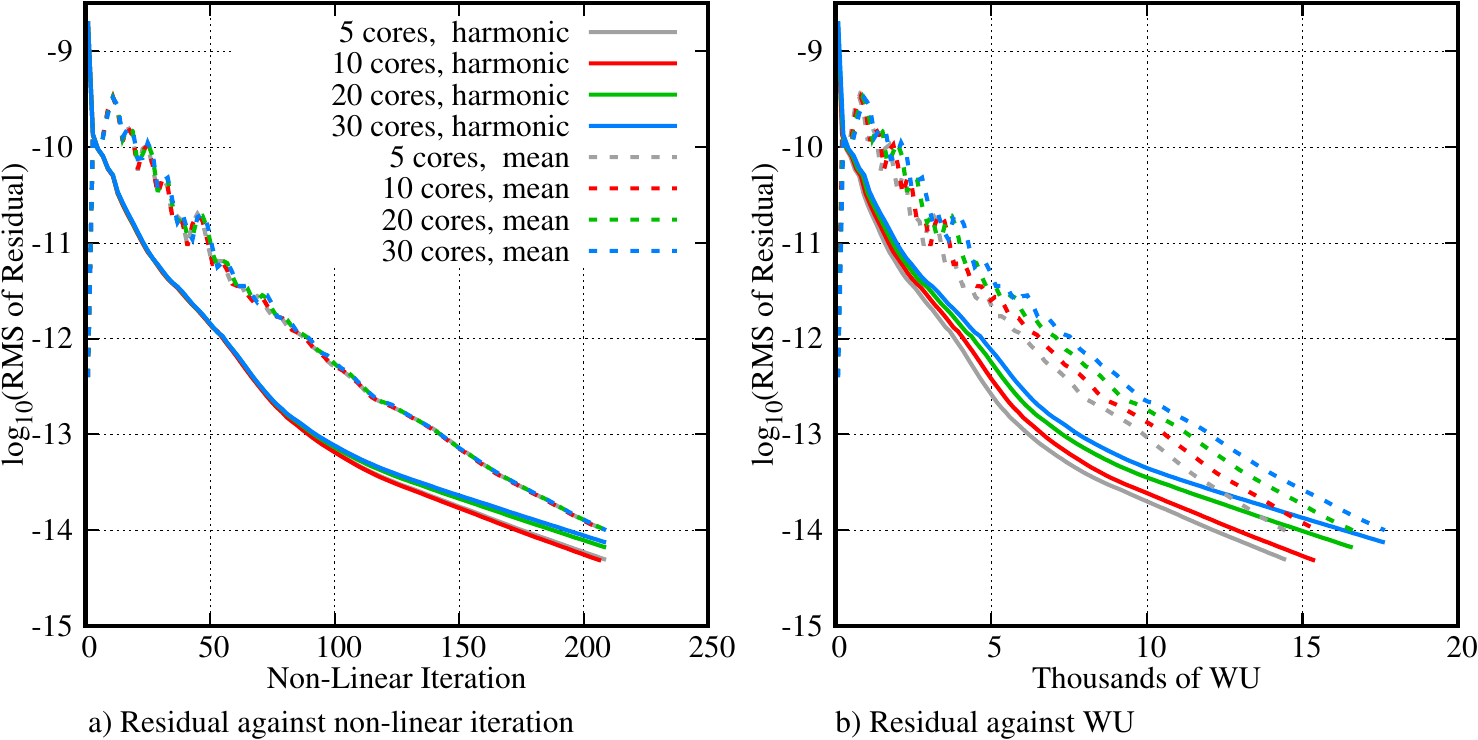}
\caption{Impact of CPU cores on the performance of the implicit scheme for the compressor case.}
\label{fig:r1_s1_r2_para}       
\end{figure}

As shown in Fig.~\ref{fig:partition_matrix}, in the parallel implementation of the implicit scheme, the entries in the Jacobian matrix that correspond to the coupling of two nodes will be ignored if they are in different partitions. Therefore, the convergence rate of the implicit scheme will depend on the number of processors. Figure~\ref{fig:r1_s1_r2_para} shows the convergence history with 5, 10, 20 and 30 CPU cores and the CFL number is set to 50. For the case with 30 CPU cores, the number of elements per core is around 1100.

In terms of nonlinear iterations, when the nonlinear residual is greater than $10^{-13}$, the number of cores only has a marginal impact on the convergence rate. When a deeper convergence is sought, adding more CPU cores leads to a slightly lower convergence rate, especially for the linearized harmonic flow. This is because the quality of the Jacobian matrix reduces as more CPUs are used. For a fixed number of linear iterations, solution updates that use more CPU cores are less accurate than those that use fewer CPU cores.

With respect to computational efficiency, which can be measured in terms of the work unit, Fig.~\ref{fig:r1_s1_r2_para} shows that adding more CPUs leads to a less efficient implicit scheme, but it has marginal effect on the number of nonlinear iterations. This is because the accuracy of the Jacobian matrix decreases as more CPU cores are used, and this will lead to a greater number of linear iterations to reach the convergence threshold of the linear residual. Hence, each nonlinear iteration becomes more expensive, but the resulting total number of nonlinear iterations can still be similar due to the fixed linear residual threshold.

\subsection*{Hotstreak Migration in A Multi-Stage Turbine}
The second case is based on the MT1 turbine stage~\cite{mt1_hotstreak} and was modified to a 1.5-stage turbine in our previous work~\cite{fnlh_multirow} to validate FNLH for hot streak migrations. The original blade count of the MT1 stage is 32 (NGV):60 (R1) and the blade counts of the redesigned 1.5-stage turbine are scaled to 30 (NGV): 60 (R1): 30 (S2). A stream tube of the geometry is simulated, and the stream tube accounts for roughly $5\%$ of the blade span. The NGV is modeled as a pair to prescribe a simplified hot streak at the NGV inlet. One layer of hexahedral elements is used to discretize the computational domain, and there are 72328 hexahedral elements in total. The average $y^{+}$ on the blade surface is around 5.  The hot streak, which spans two NGV passages, has a sinusoidal shape with an amplitude of $\Delta T_{hs}$.  The resulting inlet total temperature field $T_{\text{in}}$ is represented as the following:
\begin{equation}
    T_{\text{in}} = T_{0} + \Delta T_{hs} \sin(N_{\text{hs}} \theta)
    \label{eqn:hotstreak_def}
\end{equation}
in which $T_{0}$ is the mean total temperature at the inlet and its value is  444 K, and $N_{\text{hs}} = 15$, which means the hot streak will span two NGV passages. The measured peak-to-mean and minimum-to-mean temperature ratios are approximately 1.08 and 0.93 in the MT1 stage~\cite{mt1_hotstreak}, respectively. In this study $\Delta T_{hs} = 50$K is used. In this study, the ideal gas model is used.   For the URANS simulation, single-passage meshes are replicated in the circumferential direction to form a $12^{\circ}$ sector. 200 time steps are used for this $12^{\circ}$ sector simulation, which corresponds to 100 time steps for the stator to sweep one passage of the rotor.

The setup for the harmonics in FNLH is summarized in Table~\ref{tab::fnlh_setup_hs} and this setup is based on our previous study~\cite{fnlh_multirow}. 16 harmonics are computed in total. To model the hot streak in S2, 7 harmonics are used.  Two fundamental modes of the hot streak and their related scattered modes are computed. Their temporal frequency is related to the shaft speed $\Omega$ and its wave number is related to the linear combination of the wave number of the hot streak and the blade count of R1.

\begin{table}
\caption{FNLH setup for hot streak migration}
\label{tab::fnlh_setup_hs}
\centering
\tabcolsep=0.11cm
\begin{tabular}{llll} 
 \hline
 \hline
 Bladerow & Temporal mode & Spatial mode & Note \\ 
 \hline
  NGV & $-i\text{N}_{R1}\Omega, i\in\{1,2,...16\}$ & $i\text{N}_{R1}, i\in\{1,2,...16\}$ & exit, R1 potential field  \\ 
 \hline
 R1 & $i\text{N}_{NGV}\Omega, i\in\{1,2,...14\}$ & $i\text{N}_{NGV}, i\in\{1,2,...14\}$ & inlet, NGV wake + hot streak \\ 
      & $i\text{N}_{S2}\Omega, i\in\{1,2\}$ & $i\text{N}_{S2}, i\in\{1,2\}$ & exit, S2 potential field\\ 
       \hline
 S2 & $i\text{N}_{R1}\Omega, i\in\{1,2,...9\}$ & $i\text{N}_{R1}, i\in\{1,2,...9\}$ & inlet, R1 wake \\ 
      & $-n\text{N}_{R1}\Omega, n \in \{-1,0,1,2\}$ & $m\text{N}_{hs} + n \text{N}_{R1}, m \in \{1\}$ & inlet,Hotstreak-1 \\ 
      & $-n\text{N}_{R1}\Omega, n \in \{0,1,2\}$ & $m\text{N}_{hs} + n \text{N}_{R1}, m \in \{2\}$ & inlet,Hotstreak-2 \\ 
\hline 
\hline
\end{tabular}
\end{table}

\begin{figure}
\centering 
\includegraphics[height=7.cm]{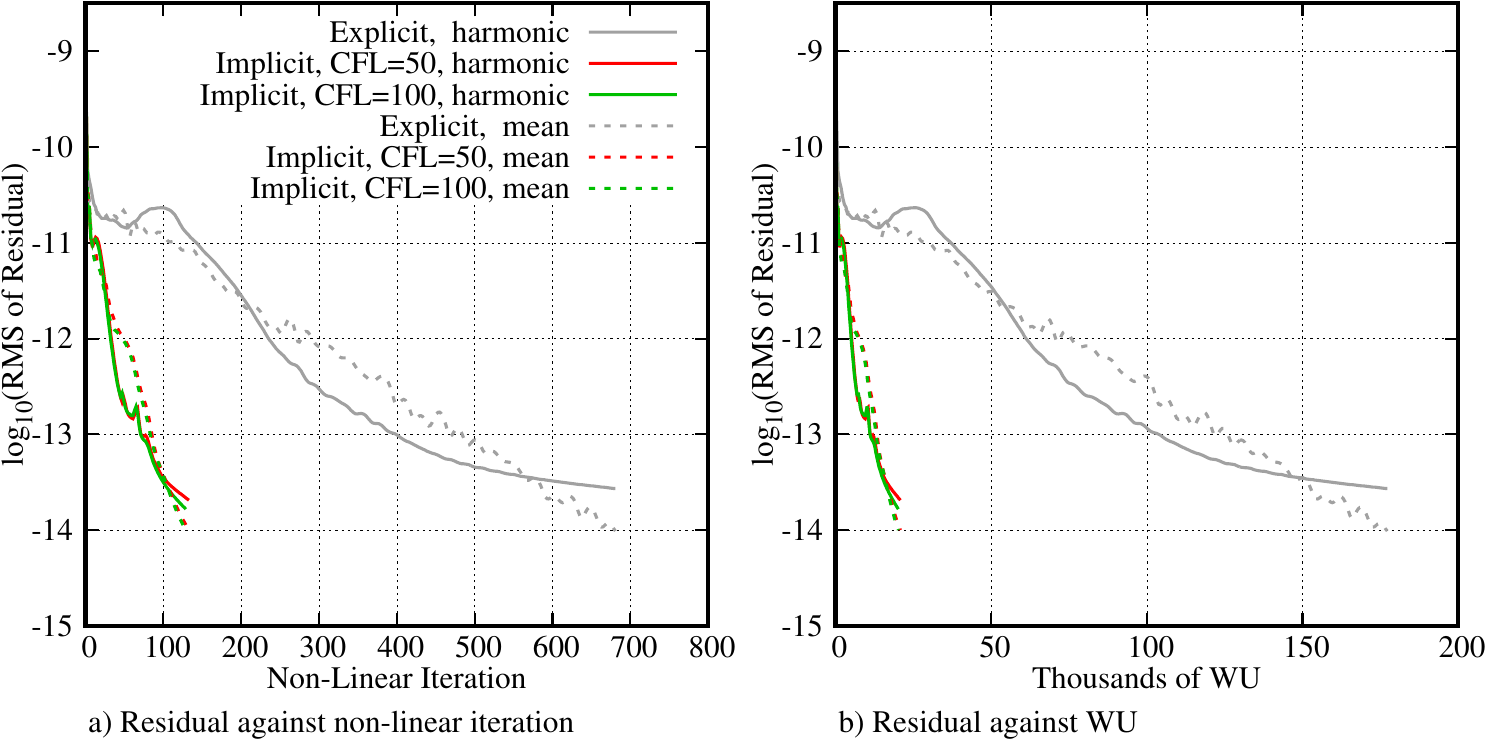}
\caption{Convergence history of explicit and implicit schemes for the turbine case.}
\label{fig:mt1_resd_rms_exp_imp}       
\end{figure}

The explicit scheme with 4 levels of MG and the implicit schemes with CFL$=50$ and CFL$=100$ are compared. The residuals of the explicit and implicit schemes are shown in Figure~\ref{fig:mt1_resd_rms_exp_imp} and 30 CPU cores are used. Similarly to the previous case, the implicit scheme shows a significant speed-up of roughly 10 compared to the explicit scheme in terms of computational time. For the explicit scheme with MG, the convergence rate is degraded as the residual decreases. This indicates that the explicit scheme with MG is not effective in damping the short-wavelength errors in the computation. However, the implicit scheme shows that it is less affected by this issue. Furthermore, for the implicit scheme, increasing the CFL from 50 to 100 only slightly improves the convergence rate. This is similar to the observation from the compressor case. 

\begin{figure}
\centering 
\includegraphics[height=7.0cm]{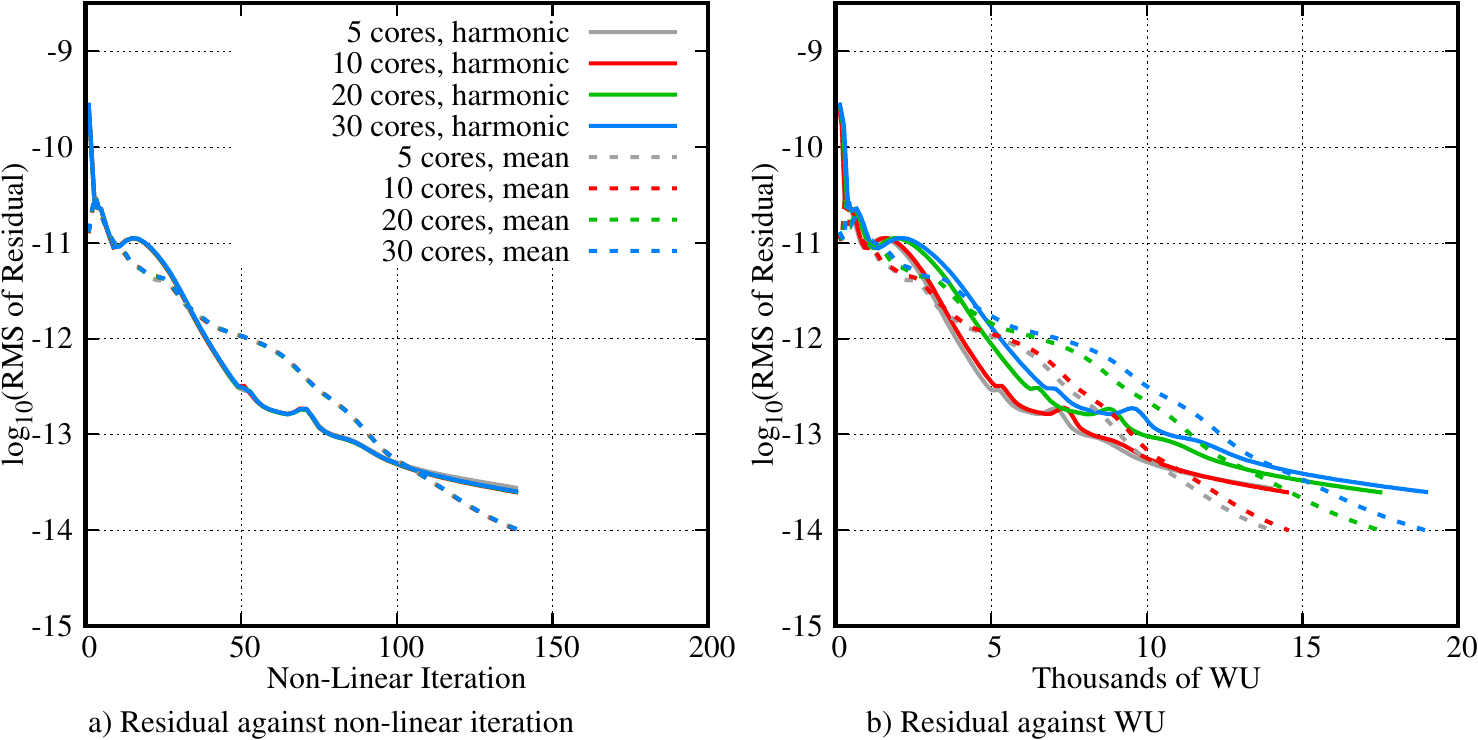}
\caption{Impact of CPU cores on the performance of the implicit scheme for the turbine case.}
\label{fig:mt1_para}       
\end{figure}

In terms of parallel performance of the implicit scheme, Fig.~\ref{fig:mt1_para} shows the convergence history with 5, 10, 20, and 30 CPU cores. For the case with 30 CPU cores, the number of elements per CPU core is roughly 2410. It can be seen that, in terms of non-linear iterations, the number of CPU cores has a negligible effect. With respect to WU, it can be seen that the implicit solver becomes slightly less efficient as the number of CPUs increases. This observation is consistent with the previous compressor case. 

\begin{figure}
\centering 
\includegraphics[height=8.cm]{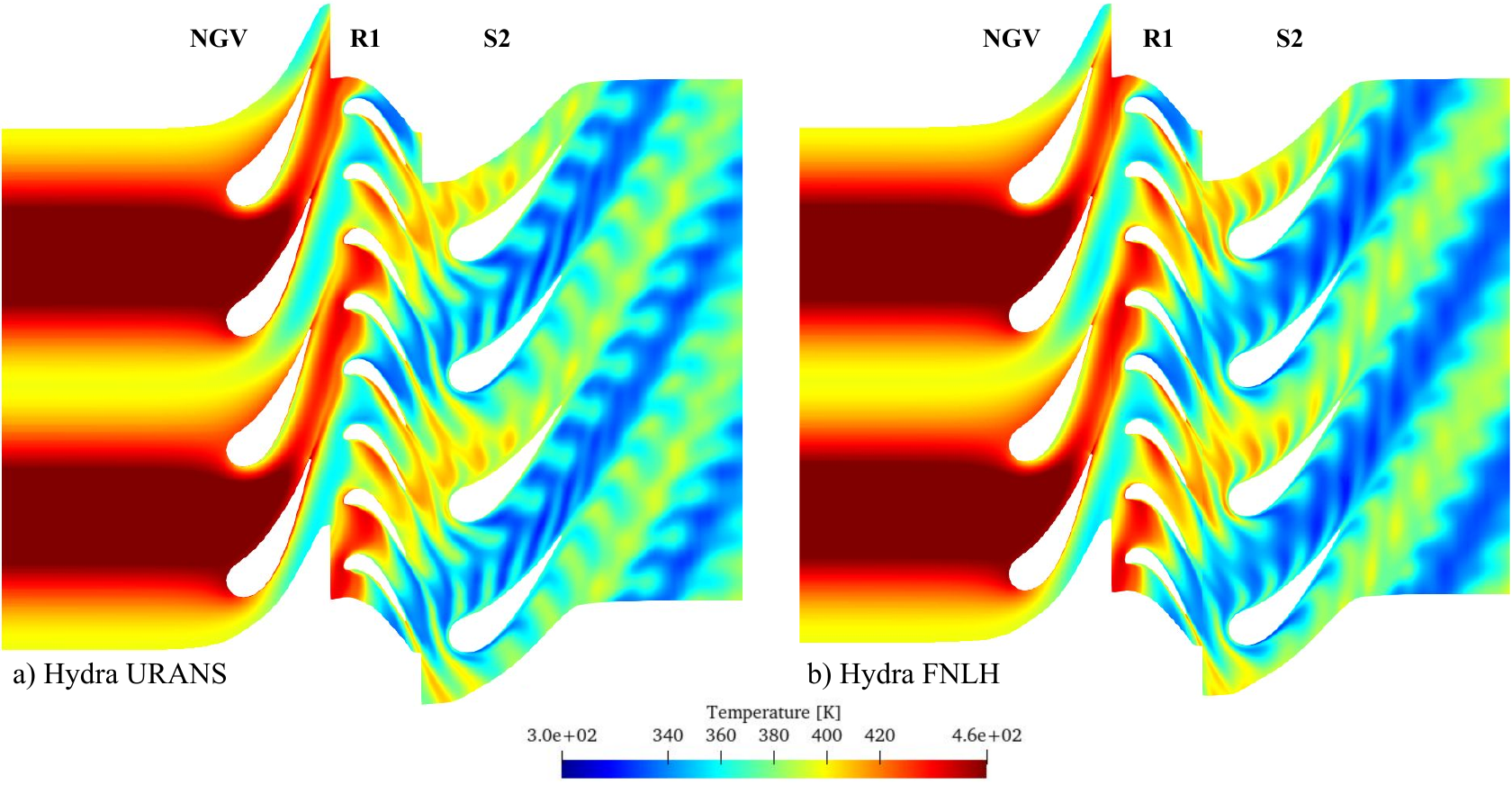}
\caption{Instataneous temperature field predicted by FNLH and URANS.}
\label{fig:mt1_mr_t}       
\end{figure}

\begin{figure}
\centering 
\includegraphics[height=8.2cm]{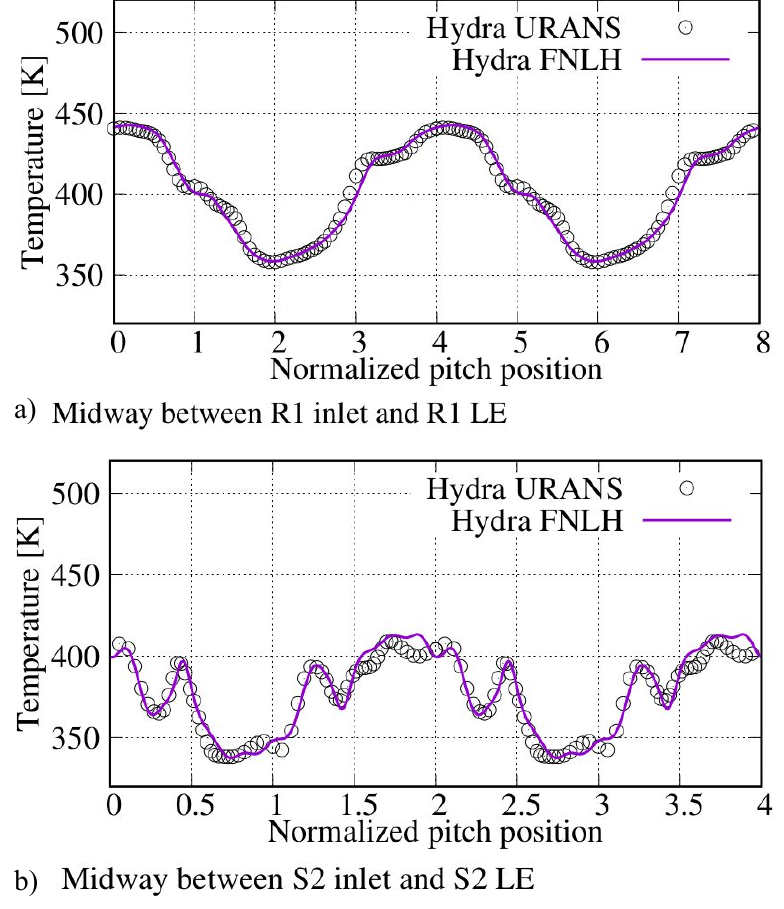}
\caption{Circumferential distribution of temperature field from FNLH.}
\label{fig:mt1_s2_le}       
\end{figure}

Figure~\ref{fig:mt1_mr_t} compares the predicted instantaneous temperature field between FNLH and URANS. For the FNLH, the implicit scheme is used. There is good agreement between URANS and FNLH on the predicted temperature field, and the hot streak is correctly transported downstream the turbine bladerows. For a quantitative comparison, Fig.~\ref{fig:mt1_s2_le} shows the circumferential distribution of the instantaneous temperature at two streamwise locations. The first is the midway between the R1 inlet and the R1 LE, and the second is the midway between the S2 inlet and the S2 LE. It can be seen that there is good agreement between the implicit FNLH and URANS. More validations of FNLH in this case can be found in previous work~\cite{fnlh_multirow}. 

\subsection*{Low Pressure Compression System of A Turbofan Engine}
\begin{figure}
\centering 
\includegraphics[height=10.5cm]{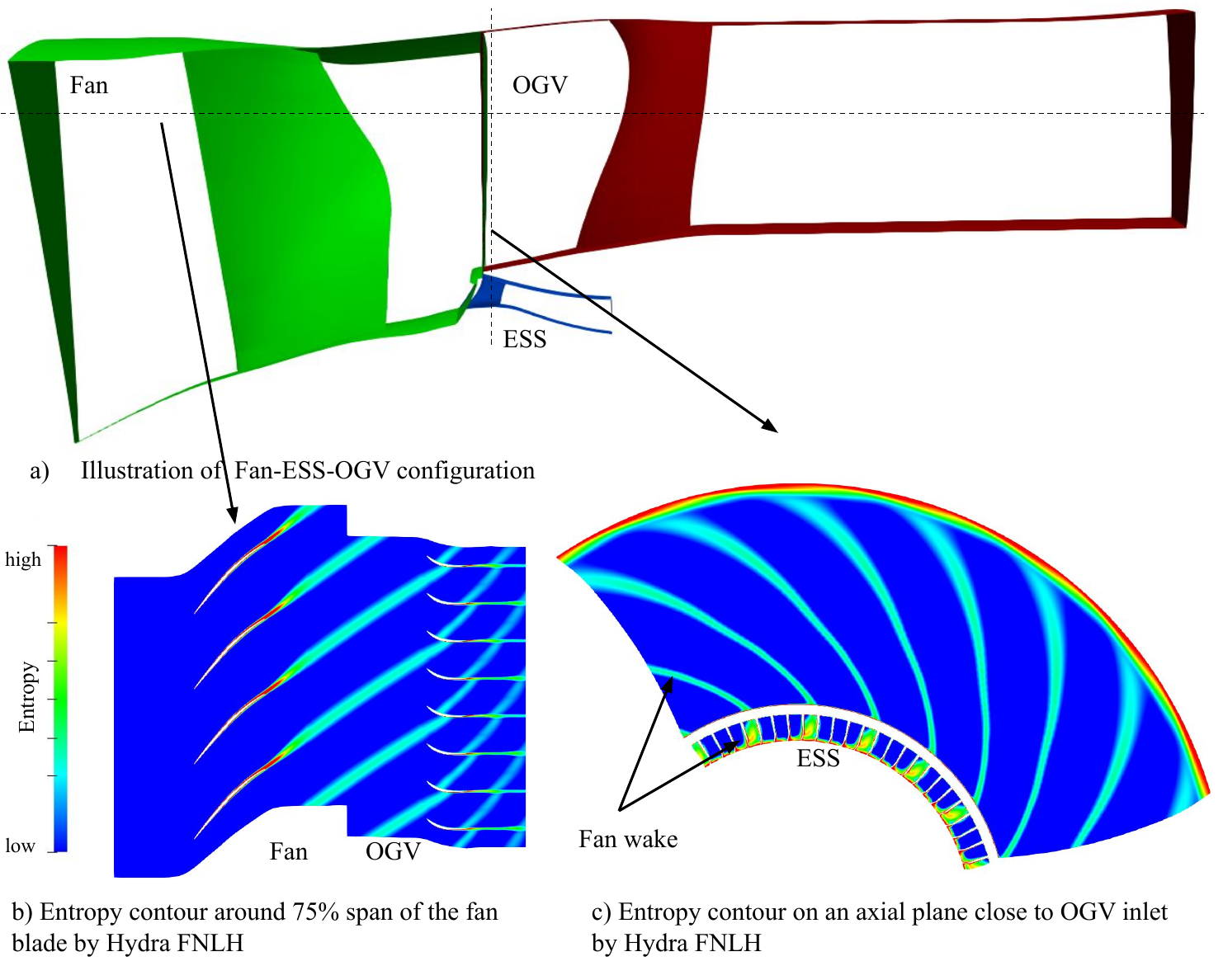}
\caption{Reconstructed instantaneous entropy of the turbo-fan LPC system.}
\label{fig:lsf4_geo}       
\end{figure}

The third test case is the low pressure compression (LPC) system of a high bypass ratio turbofan engine. It consists of the fan, engine section strut (ESS), and the outlet guide vane (OGV). The computational domain of this case is shown in Fig.~\ref{fig:lsf4_geo}. The exit boundary of the fan blade is split into two parts, one is connected to ESS and the other one is linked to OGV. The flow in the fan passage features the passage shock towards the casing and strong three-dimensional secondary flows.  The fan blade is discretized with 110 layers of hexahedral in the spanwise direction and has 2077236 hexahedral cells in total.   The FNLH setup for this case is summarized in Table~\ref{tab::fnlh_setup_lsf4}. 8 harmonics are computed in total. In particular, for the fan blade, 4 harmonics are used to represent the potential field from ESS and OGV, respectively. For ESS and OGV, 8 harmonics are used to model the fan wake. Figure~\ref{fig:lsf4_geo} shows the instantaneous entropy reconstructed around the $75\%$ span of the fan and it can be seen that the fan wake is transmitted into the computational domain of OGV without attenuation on the fan-OGV interface. Figure~\ref{fig:lsf4_geo} also shows the instantaneous entropy distribution in an axial plane located close to the OGV inlet. It can be seen that the fan wake is correctly reconstructed in the computational domain of OGV and ESS. As the focus of the paper is on the computational performance of the generalized RK scheme, more detailed quantitative validations of the URANS will be presented in future work. 

\begin{table}
\caption{FNLH setup for the low pressure compression system of a turbofan}
\label{tab::fnlh_setup_lsf4}
\centering
\tabcolsep=0.11cm
\begin{tabular}{cccc} 
 \hline
 \hline
 Bladerow & Temporal mode & Spatial mode & Note \\ 
 \hline
  Fan & $-i\text{N}_{ESS}\Omega_{Fan}, i\in\{1,2,3,4\}$ & $i\text{N}_{ESS}, i\in\{1,2,3,4\}$ & exit, ESS potential field  \\ 
      & $-i\text{N}_{ESS}\Omega_{Fan}, i\in\{1,2,3,4\}$ & $i\text{N}_{OGV}, i\in\{1,2,3,4\}$ & exit, OGV potential field  \\ 
 \hline
 ESS     & $i\text{N}_{Fan}\Omega_{Fan}, i\in\{1,2,...8\}$ & $i\text{N}_{Fan}, i\in\{1,2,...8\}$ & inlet, Fan wake\\ 
       \hline
 OGV & $i\text{N}_{Fan}\Omega_{Fan}, i\in\{1,2,...8\}$ & $i\text{N}_{Fan}, i\in\{1,2,...8\}$ & inlet, Fan wake \\ 
 
\hline 
\hline
\end{tabular}
\end{table}

\begin{figure}
\centering 
\includegraphics[height=6.cm]{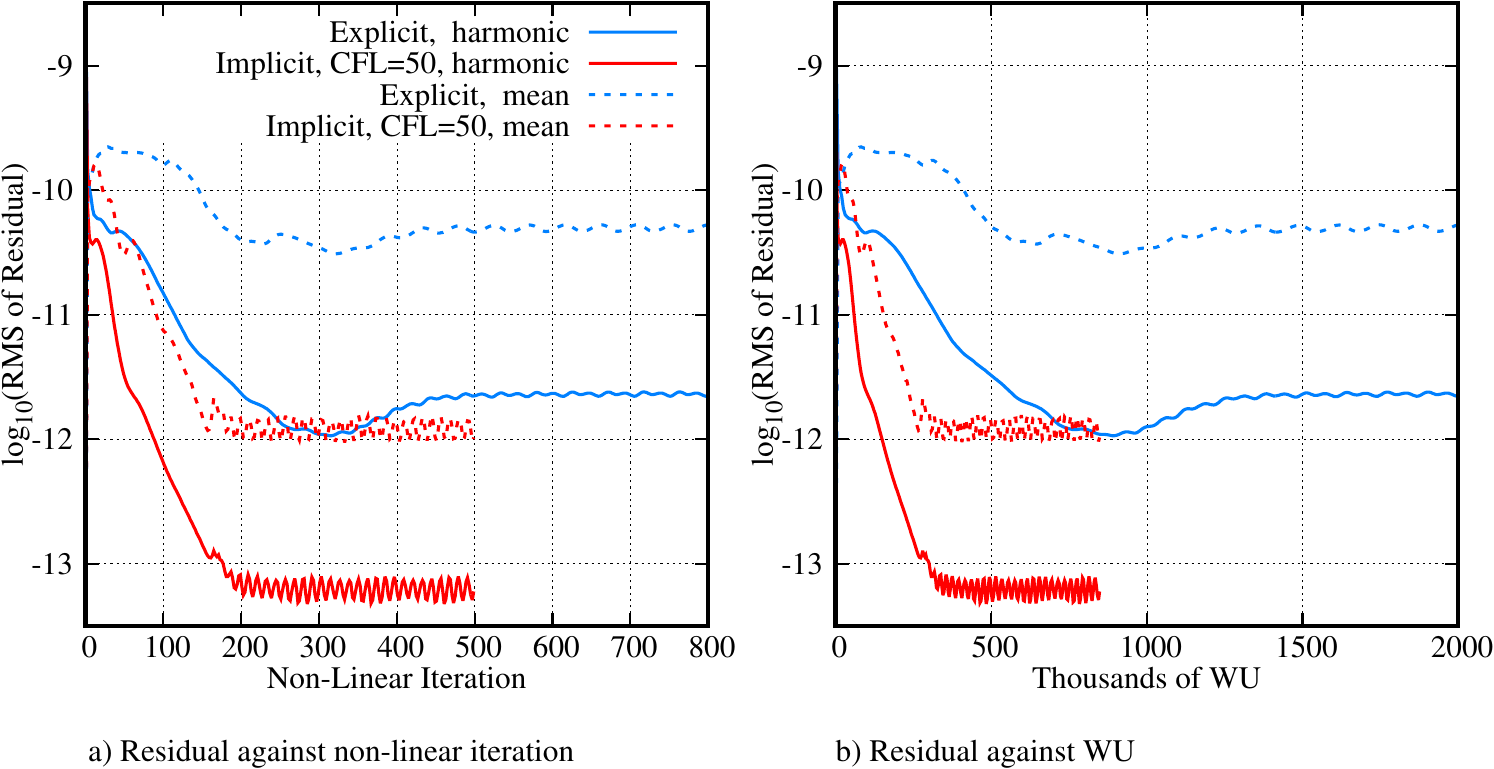}
\caption{Comparisons of explicit and implicit schemes for the LPC system. }
\label{fig:lsf4_rms_exp_imp}       
\end{figure}

Figure~\ref{fig:lsf4_rms_exp_imp} shows the comparisons of the convergence history of explicit and implicit RK schemes. 40 CPU cores are used. The CFL number of the implicit scheme is 50 and 4 levels of MG and V-cycle are used by the explicit scheme. It can be seen that the implicit scheme shows roughly 8x speed-up compared to the explicit scheme in terms of computational time. This is consistent with the previous two cases. In addition, the residuals of the explicit schemes stagnate to a higher level, while the implicit scheme shows a better convergence rate and can reach a deeper convergence compared to the explicit scheme.  This shows that for challenging industrial cases, the implicit scheme not only shows a faster convergence rate, but also reaches a deeper convergence compared to the explicit scheme. 

\subsection*{Remarks on Computational Resources}
Table~\ref{tab::memory_consumption} compares the memory consumption of the explicit and implicit RK schemes for the test cases. It can be seen that the memory consumption of the implicit solver is roughly half of the explicit solver using 4 MG levels. Therefore, the implicit RK scheme not only delivers faster and deeper convergence compared to the explicit scheme, but also has significantly less memory consumption. The memory consumption of the corresponding whole annulus URANS simulations is also summarized in Table~\ref{tab::memory_consumption}. The memory consumption of FNLH is roughly two orders of magnitude smaller compared to the whole annulus URANS.

\begin{table}
\caption{Memory consumption ( in GB ) of the test cases}
\label{tab::memory_consumption}
\centering
\tabcolsep=0.11cm
\begin{tabular}{cccc} 
 \hline
 \hline
         & Explicit (4-MG) & Implicit & URANS ($360^{\circ}$)  \\ 
 \hline
  Case 1 & 3.9  & 2.1  & 309.6   \\ 
 \hline
 Case 2  & 7.4  & 3.2  & 378.1  \\ 
       \hline
 Case 3 & 56.5  & 32.6  & 960.2   \\ 
 
\hline 
\hline
\end{tabular}
\end{table}

\begin{figure}
\centering 
\includegraphics[height=5.25cm]{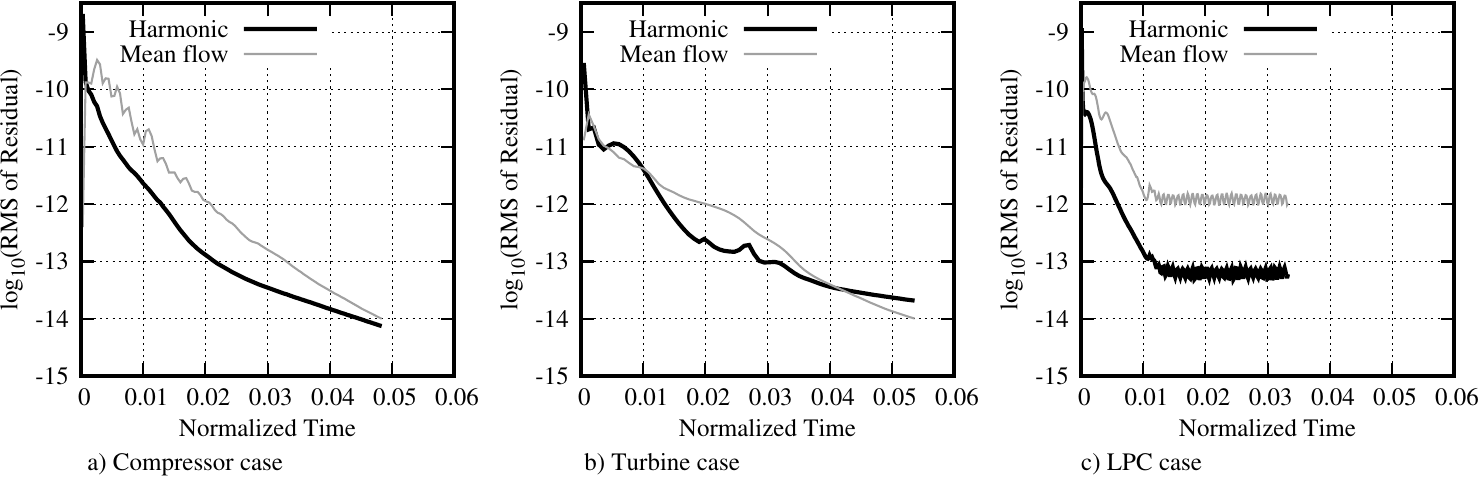}
\caption{Computational efficiency of FNLH compared to URANS.}
\label{fig:normalized_time}       
\end{figure}

Figure~\ref{fig:normalized_time} shows the convergence history of the implicit scheme (CFL = 50) for the test cases, but the x-axis of the plots is the normalized time, which is defined as the wall clock time of the FNLH simulation divided by the wall clock time of the execution of the full annulus URANS simulation for one revolution. It is typical to run 4-5 revolutions before collecting data from URANS. Figure~\ref{fig:normalized_time} shows that the implicit FNLH is at least two orders of magnitude more efficient than the corresponding full-annulus URANS simulations in all test cases.  

\section*{Conclusions and Future Work}
A unified formulation has been developed and implemented in the industrial CFD solver Hydra to accelerate the convergence of the FNLH method. A comparative study of the resulting explicit and implicit schemes has been undertaken. The proposed method is memory efficient as memory consumption depends weakly on the number of harmonics to compute. For the implicit scheme, the complex-number representation of the matrix is used for the matrix factorization for the purpose of reusing the existing lookup table for that of the mean flow. This significantly reduces the effort of implementing the implicit FNLH scheme.  


The implicit FNLH scheme is found to consistently show a speedup of 7-10 compared to the explicit schemes with multigrid and this speed-up is comparable to the value reported by Misev~\cite{misev_thesis} for the baseline steady flow solver. This demonstrates the optimal performance of the current implementation of the implicit FNLH method.  Furthermore, the memory consumption of the implicit scheme consumes approximately $50\%$ of that of the explicit scheme with 4 levels of multigrid. The parallel performance of the implicit scheme is also studied and shows that the implicit solver is slightly less efficient as more CPU cores are used. 



Future work includes removing the restriction on frozen turbulence and considering the cross-coupling of the harmonics. These treatments could have an impact on the convergence and robustness of the FNLH method. Besides, the robustness of the proposed procedure for cases with strong flow separation is not fully investigated in this paper. As the flow approaches the stall condition, it can generate unsteady flow disturbances whose frequencies are not relevant to the shaft speed. As harmonic methods require the specification of temporal and spatial modes a-priori, this poses further challenges to harmonic methods. Therefore, the implicit solver needs to be enhanced to overcome this challenge,  and this will be investigated in future work.

\section*{Acknowledgments}
The authors are grateful to Rolls-Royce plc and Aerospace Technology Institute for funding this work under the CORDITE project and granting permission for its publication. The authors would like to thank Dr. David Liliedahl from Rolls-Royce Corporation (USA) for discussions on the implicit solver of Hydra. The authors would like to thank Dr. Paolo Adami for his support of this work.

\begin{nomenclature}

\entry{\textbf{Latin symbols}}{}
\entry{$a$}{speed of sound}
\entry{$I$}{imaginary unit, $\sqrt{-1}$, unit matrix}
\entry{$l$}{circumferential wave number, harmonic index}
\entry{$m$}{circumferential wave number}
\entry{$n$}{circumferential wave number}
\entry{$N$}{number of passages}
\entry{$N_h$}{number of harmonics}
\entry{$u$}{velocity}
\entry{$U$}{conservative flow variable}
\entry{$p$}{pressure}
\entry{$W$}{primitive flow variable}
\entry{$\mathbf{R}$}{residual}
\entry{$s$}{entropy}
\entry{$t$}{time}
\entry{$T$}{temperature}

\entry{\textbf{Greek symbols}}{}
\entry{$\rho$}{density}
\entry{$\omega$}{angular frequency}
\entry{$\Omega$}{shaft rotational speed}

\entry{\textbf{Subscripts, superscripts and operators}}{}
\entry{$i$}{spatial index}
\entry{$j$}{spatial index}
\entry{$k$}{space navigation index}
\entry{$l$}{Fourier index}
\entry{$'$}{fluctuation relative to the passage-averaged quantity}
\entry{$''$}{fluctuation relative to the Favre-averaged quantity}
\entry{$\bar{}$}{time-average}
\entry{$\tilde{}$}{Favre-average}
\entry{$\hat{}$}{Fourier coefficient}

\entry{\textbf{Acronym}}{}
\entry{CFD}{Computational Fluid Dynamics}
\entry{DF}{Deterministic Flux}
\entry{DOF}{Degree of freedom}
\entry{FNLH}{Favre-averaged Non-Linear Harmonic}
\entry{HB}{Harmonic Balance}
\entry{HS}{Hot streak}
\entry{IBPA}{Inter-Blade Phase Angle}
\entry{MG}{Multi-Grid}
\entry{NGV}{Nozzel Guide Vane}
\entry{NLH}{Non-Linear Harmonic}
\entry{NS}{Navier-Stokes}
\entry{RANS}{Reynolds-Averaged Navier-Stokes}
\entry{RMS}{Root-mean-square}
\entry{URANS}{Unsteady Reynolds-Averaged Navier-Stokes}
\entry{WU}{Work unit}
\end{nomenclature}

\bibliographystyle{asmems4}

\bibliography{asme2ej}


\end{document}